\documentclass[10pt]{article}

\usepackage[utf8]{inputenc}   
\usepackage[T1]{fontenc}      
\usepackage[french]{babel}    
\usepackage{amsmath, amssymb} 
\usepackage{graphicx}         
\usepackage{geometry}         
\usepackage{subcaption}
\usepackage{amsthm}  
\usepackage{tikz}
\usepackage{amsmath}
\usepackage{pgfplots}
\usetikzlibrary{matrix} 
\usetikzlibrary{arrows.meta, decorations.markings}
\usepackage{blkarray} 
\usepackage{float}
\usepackage{amsmath,amssymb}
\usepackage{mathtools} 

\usepackage{subcaption} 
\usepackage{amsmath,amssymb,amsthm}

\theoremstyle{plain}
\newtheorem{proposition}{Proposition}[section]
\newtheorem{lemma}[proposition]{Lemma}
\newtheorem{convention}{Convention}

\theoremstyle{remark}
\newtheorem{remark}[proposition]{Remark}

\usepackage[hidelinks]{hyperref}
\usepackage[nameinlink,noabbrev]{cleveref} 

\usepackage[symbol]{footmisc}

\usetikzlibrary{patterns}

\makeatletter
\newtheoremstyle{boldthm}  
  {3pt}{3pt}               
  {\itshape}               
  {}                       
  {\bfseries}              
  {.}                      
  {0.5em}                  
  {}                       
\makeatother

\theoremstyle{boldthm}
\newtheorem{theorem}{Theorem}

\title{The eigenvector bead process}
\author{
  Antonin Barbe,
  Benjamin De Bruyne,
  Romain Allez\footnotemark[2] 
}

\date{\today}

\begin{document}

\maketitle
\footnotetext[2]{Qube Research \& Technologies, Paris, France. \texttt{antonin.barbe@qube-rt.com},\texttt{benjamin.debruyne@qube-rt.com},\\ \texttt{romain.allez@qube-rt.com}}

\begin{abstract}
We investigate the overlap matrix between the eigenvectors of a Wigner matrix $H_{N+K}$ of size $(N+K)\times(N+K)$ and those of its principal minor $H_N$ of size $N\times N$, for both the real symmetric ($\beta=1$) and complex Hermitian ($\beta=2$) ensembles, in the regime where $N \to \infty$ while $K$ remains fixed. 

Our analysis yields two main results. 
(i) In the \emph{bulk} of the spectrum, an eigenvector of $H_{N+K}$ associated with an eigenvalue at energy level $E$ projects primarily onto eigenvectors of $H_N$ located at the same local spectral level. This phenomenon, which we call \emph{local projection}, highlights a robust stability of the eigenbasis under matrix growth. 
(ii) At the \emph{spectral edge}, the change of basis between the leading eigenspaces of consecutive minors is asymptotically governed by a random antisymmetric perturbation of order $N^{-1/3}$. In both cases, we provide the asymptotic law of the overlaps expressed in terms of the Airy and Sine kernels.

We further extend our analysis to the case of Wishart matrices, that is, sample covariance matrices of the form $W = X^{\!\top} X$, where $X \in \mathbb{R}^{T \times N}$ is a matrix with i.i.d.\ random entries.
We establish analogous results for the overlaps between eigenvectors of consecutive minors of $W$, both in the bulk and at the spectral edges (soft and hard). 
The limiting laws share the same universal structure as in the Wigner case, up to explicit constants depending on the aspect ratio $q = N/T$. 
This demonstrates the universality of the eigenvector overlap process across distinct random matrix ensembles.
\end{abstract}

\section{Introduction}

Random matrix theory has uncovered a rich array of universal behaviors in the spectral properties of large symmetric matrices, commonly known as Wigner matrices \cite{wigner1955,wigner1957}. While the global eigenvalue distribution (given by Wigner’s semicircle law) \cite{mehta2004} and local eigenvalue statistics (e.g. Wigner–Dyson bulk spacings and Tracy–Widom edge fluctuations) have been thoroughly studied and shown to be largely universal across matrix ensembles \cite{tracy1994,tracy1996,dyson1962,soshnikov1999,johansson2001,taovu2010,gorin-shkolnikov2015,majumdar-schehr2014}, the structure and evolution of eigenvectors remain far less understood \cite{erdos2009,allez-bouchaud2014,bun-bouchaud-potters2018}.

Let \(H_N\) denote an \(N\times N\) real symmetric Wigner (GOE) matrix, and for \(1\le n\le N\) let \(H_n\) be its principal \(n\times n\) minor (taken on the first \(n\) coordinates). In the macroscopic regime—namely, when \(n,N\to\infty\) with \(n/N\to p\in(0,1)\)—the mean squared overlaps between an eigenvector of $H_n$ and one of $H_N$ scale as $O(1/N)$ as $N\to \infty$. The eigenvectors of $H_n$ are spread across the whole spectrum of $H_N$ with a density profile that depends on the spectral locations of the eigenvalues associated to the eigenvectors \cite{attal2025,attal-allez2025}.

The microscopic scaling limit---where \(n=N-K\) with fixed \(K\ge 1\)---remains an open problem. In this work we address this problem by studying how the eigenvectors of a Wigner matrix evolve as the matrix grows, one row and one column at a time. This natural enlargement defines a discrete Markovian dynamic, since each such rank-one update perturbs the matrix and thereby reshapes its eigenvalues and eigenvectors \cite{bunch1978original}. In contrast to the macroscopic limit, an eigenvector of $H_n$ is expected to be localized onto the eigenvectors of $H_N$ associated to eigenvalues at the same energy level. In addition, while the overlaps in the macroscopic case are self-averaging, the ones in the microscopic case remain $O(1)$ random variables in the limit of $N\to\infty$.

For the integrable Gaussian case, a great deal is known about the eigenvalue minor process. In fact, for the Gaussian Unitary Ensemble (GUE), the joint distribution of eigenvalues of all principal minors (the GUE corner process) is exactly solvable. Baryshnikov \cite{baryshnikov2001} first observed a connection between the GUE minors and last passage percolation models, and later Johansson and Nordenstam \cite{johansson-nordenstam2006} showed that the GUE corner process is a determinantal point process with explicit correlation kernels, converging in the bulk to the so-called “bead” process \cite{boutillier2009,najnudel-virag2021}. In general for Wigner matrices (with non-Gaussian entries), such exact formulas are unavailable and one must instead rely on universality methods \cite{huang2022}. Huang’s recent work, entitled \emph{Eigenvalues for the Minors of Wigner Matrices}, demonstrated that the eigenvalues of Wigner minors indeed form a well-behaved Markov chain as the size grows \cite{huang2022}. At the spectral edge, in particular, the adjacent gaps in this Wigner corner process behave like a random walk with independent Gamma-distributed steps \cite{huang2022,gorin-shkolnikov2017}. This remarkable “edge decoupling” phenomenon was predicted in the $\beta$-ensemble setting (Hermite $\beta$ corner process) by Gorin and Shkolnikov, who proved that the spacings between extreme eigenvalues of adjacent minors converge to i.i.d. gamma variables \cite{gorin-shkolnikov2017}.

Further structural insights have emerged about the transition between fully correlated and fully decorrelated regimes in the minor process. Forrester and Nagao \cite{forrester2011} provided an explicit description of the \emph{critical} regime in the GUE minor process—namely when the size difference \(K:=N-n\) between minors satisfies \(K=\Theta\!\left(N^{2/3}\right)\)—giving a closed-form correlation kernel for the resulting interpolating determinantal point process. Very recently, Bao et al. \cite{bao2025} proved that this decorrelation transition is universal for Wigner matrices as well: even for non-Gaussian Wigner ensembles, the correlation between the eigenvalues of a full $N\times N$ matrix and those of a smaller $n\times n$ principal minor undergoes the same critical decay as $N-n$ grows, in agreement with the integrable GUE predictions. Notably, this result also extends to the overlaps of eigenvectors, showing that the eigenvectors of a smaller Wigner minor become asymptotically orthogonal to those of the full matrix beyond the critical minor size \cite{bao2025}. We remark that related questions for dynamically growing matrices have been studied in the context of Dyson Brownian motion; see \cite{sodin2015} for results on the time-dependent Wigner corner process.

\bigskip
Let us briefly recall the following result on the eigenvalue process due to Huang \cite{huang2022}, on which the results of this paper rely on.  
Fix integers $K, \ell \geq 1$. Let $H_{N+K}$ be a Wigner matrix of size \((N+K)\times (N+K) \), and let 
$\lambda_1^{(s)}\leq\ldots\leq\lambda_{N+s}^{(s)}$ denote the eigenvalues of the principal minor of size $(N+s)\times(N+s)$. 
The array
\[
\left(
N^{1/6}\bigl(\lambda_1^{(0)} + 2\sqrt{N}\bigr),\quad 
N^{1/6}\bigl(\lambda_2^{(0)} + 2\sqrt{N}\bigr),\quad 
\ldots,\quad 
N^{1/6}\bigl(\lambda_\ell^{(0)} + 2\sqrt{N}\bigr)
\right)
\]
converges in distribution, as $N \to \infty$, to the Tracy--Widom$_\beta$ distribution.  
\noindent
Moreover, the array
\[
\left( 
\sqrt{N}\bigl(\lambda_1^{(s)} - \lambda_1^{(s+1)}\bigr),\quad 
\sqrt{N}\bigl(\lambda_2^{(s)} - \lambda_2^{(s+1)}\bigr),\quad 
\ldots,\quad 
\sqrt{N}\bigl(\lambda_\ell^{(s)} - \lambda_\ell^{(s+1)}\bigr)
\right)_{0 \leq s < K}
\]
converges in distribution, as $N \to \infty$, to a random array with independent Gamma-distributed entries, with density 
\[
C_\beta \, x^{\beta/2 - 1} e^{-\beta x/2},
\quad \text{where } C_\beta = \frac{(\beta/2)^{\beta/2}}{\Gamma(\beta/2)}.
\]
In addition, the two arrays are asymptotically independent.  

\medskip
In the bulk of the spectrum, at an energy level $E \in (-2,2)$, the following point process defined by the eigenvalues in the neighbourhood of the level $E$:
\[
\Theta^{(s)}_N = \sum_{j=1}^{N+s} \delta_{\sqrt{N(4 - E^2)} \, (\lambda_j^{(s)} - E\sqrt{N})}, 
\quad 0 \leq s \leq K,
\]
converges in distribution as $N \to \infty$ to a Markov chain $\Theta^{(s+1)} = \mathcal{D}(\Theta^{(s)}, \Gamma^{(s)}, h)$ given in (1.4) in \cite{huang2022}.

%
%

\medskip
In the present work, we build on the key ideas of Huang \cite{huang2022} 
to describe the behavior of the eigenvectors overlaps in the microscopic regime. In the \emph{bulk}, we show that the eigenvectors remain spectrally aligned under matrix growth: 
an eigenvector at a given spectral location projects mainly onto eigenvectors from 
the previous minor near the same position. 
At the \emph{edge}, the regime is different: as the spectral density vanishes, 
leading eigenvectors of successive minors undergo small random antisymmetric 
rotations, close to the identity but occurring at scale $N^{-1/3}$.

\section{Main Results}

In this section, we present the main results of the paper. We first analyze the evolution of eigenvector overlaps under matrix growth in the Wigner ensemble, following the theoretical framework introduced in \cite{huang2022}. We then extend the analysis to the Wishart ensemble, showing that an analogous set of master equations governs the dynamics in this case as well. Despite their structural differences, both ensembles display a similar universal behavior of eigenvectors under infinitesimal matrix enlargement. The results for the Wigner and Wishart cases are presented successively.

\subsection{Wigner Ensemble}

\medskip
\noindent
\textbf{Wigner matrices.}
\label{par:wigner-matrices}
A Wigner matrix $H_N = [h_{ij}]_{1\le i,j\le N}$ is a real symmetric or complex Hermitian $N\times N$ random matrix such that the upper–triangular entries $h_{ij}$, for $i\le j$, are independent, have mean zero, and satisfy $\mathbb{E}\!\left[|h_{ij}|^2\right]=1$. In the Hermitian case, we assume $\mathrm{Re}\,h_{ij}$ and $\mathrm{Im}\,h_{ij}$ are independent and, for $i<j$, 
$\mathbb{E}\!\left[(\mathrm{Re}\,h_{ij})^2\right]=\mathbb{E}\!\left[(\mathrm{Im}\,h_{ij})^2\right]=\tfrac12$. We also require a finite fourth moment:
\[
\mathbb{E}\!\left[|h_{ij}|^4\right] < \infty.
\]

\bigskip
\noindent
Let $H_N\in\mathcal S_N(\mathbb R)$ be real symmetric with simple spectrum
$\lambda^{(N)}_1<\dots<\lambda^{(N)}_N$ and orthonormal eigenbasis
$\{\textbf{u}^{(N)}_j\}_{j=1}^N$.
Given a vector $\mathbf g=(g_1,\dots,g_N)^\top$ expressed in this basis, and a scalar $g_{N+1}\in\mathbb R$,
consider the rank-one extension
\[
H_{N+1}\;=\;
\begin{pmatrix}
H_N & \mathbf g\\
\mathbf g^\top & g_{N+1}
\end{pmatrix}\in \mathcal S_{N+1}(\mathbb R),
\qquad
g_i^{(N)}:=\langle \mathbf u_i^{(N)},\mathbf g\rangle\ (1\le i\le N).
\]
For a unit eigenvector $\mathbf u^{(N+1)}_i$ of $H_{N+1}$ associated with $\lambda^{(N+1)}_i$, define the overlaps \(\Omega^{(N)}_{i,j} \) and the eigenvalues gaps \( \Delta_{i,j}^{(N)} \)
\[
\Omega^{(N)}_{i,j}\;:=\;\big\langle \mathbf{u}^{(N+1)}_i,\,\mathbf{u}^{(N)}_j\big\rangle,
\qquad
\Omega^{(N)}_{i,N+1}:=\big(\mathbf{u}^{(N+1)}_i\big)_{N+1},
\]
\[
\Delta_{i,j}^{(N)}:=\lambda^{(N+1)}_i-\lambda^{(N)}_j \qquad (1\le i\le N+1,\ 1\le j\le N).
\]

%
\noindent
The Stieltjes transform of the discrete measure supported at the points
\((\mu_j)_{j\in\mathbb{Z}}\) with weights \(|g_j|^2\) is defined, for
\(z\in\mathbb{C}\setminus\{\mu_j:j\in\mathbb{Z}\}\), by
\begin{equation}\label{eq:secular-S}
S_{\mu,g^2}(z) := \sum_{j\in\mathbb{Z}} \frac{|g_j|^2}{\mu_j - z}.
\end{equation}

\begin{convention}[Inverse branches of the Stieltjes transform]\label{conv:inverse-stieltjes-branches}
Let \((\mu_j)_{j\in\mathbb{Z}}\) be strictly increasing and let \(|g_j|^2>0\).
For \(z\notin\{\mu_j:j\in\mathbb{Z}\}\), define
\(
S_{\mu,g^2}(z)
\)
For any \(h\in\mathbb{R}\), the equation \(S_{\mu,g^2}(z)=h\) admits a unique solution
on each interval \(I_j:=(\mu_j,\mu_{j+1})\), denoted by \(S_{\mu,g^2}^{-1}(h)_j\).
We write \((S_{\mu,g^2}^{-1})'\) for the derivative of the inverse map with respect to \(h\),
and \((S_{\mu,g^2}^{-1})'_i\) for its value on the \(i\)-th branch.
\end{convention}

\noindent
We now present our main results on the limiting law of the overlap matrix $\Omega^{(N)}$ both in the edge and in the bulk of the spectrum.  


\bigskip

\begin{theorem}[Edge case]\label{thm:edge-case}
Fix an integer \(\ell\ge 1\).
For \(1 \leq i,j \leq \ell \), the overlap matrix converges to the identity matrix, i.e.~$\lim_{N\to\infty}(\Omega_{i,j}^{(N)})_{1\leq i,j\leq \ell}=I_{\ell}$ (in agreement with \cite{bao2025}). Furthermore, the leading-order correction scales as $O(N^{-1/3})$ and is given by:
\begin{equation}\label{eq:edge-edge}
  N^{1/3}\bigl(\Omega^{(N)}_{i,j}-(I_{\ell})_{i,j}\bigr)\;
  \mathrel{\underset{N\to\infty}{\overset{d}{\sim}}}\;
  \begin{cases}
    \dfrac{\overline{g_i}\, g_j}{\alpha_j-\alpha_i}, & i\neq j,\\[0.6ex]
    0, & i=j,
  \end{cases}
  \qquad 1\leq i,j\leq \ell.
\end{equation}

\noindent
Here \((\alpha_i)_{i=1}^\ell\) follow the Airy$_{\beta}$ point process, and \((g_i)_{i=1}^\ell\) are i.i.d.\ standard Gaussians \(\mathcal{N}_\beta(0,1)\) independent of \((\alpha_i)\).
\end{theorem}

\begin{theorem}[Bulk case]\label{thm:bulk-case}
We recall Convention~\ref{conv:inverse-stieltjes-branches}. Fix \(E\in(-2,2)\) and define the anchor \\
\(
i_N(E)\;:=\;\min\bigl\{\,k\in\{1,\dots,N\}\ :\ \lambda^{(N)}_{k}\ \ge E\sqrt{N}\,\bigr\}.
\) Then the overlaps satisfy
\begin{equation}\label{eq:bulk-near-anchored}
\Omega^{(N)}_{\,i_N(E)+u,\;i_N(E)+v}
\;\mathrel{\underset{N\to\infty}{\overset{d}{\sim}}}\;
\frac{g_v}{\,S_{\mu,\,g^2}^{-1}(h)_u-\mu_v\,}\,
\sqrt{\bigl(S_{\mu,\,g^2}^{-1}\bigr)'(h)_u},
\qquad h=-\frac{E}{2\sqrt{4-E^2}}.
\end{equation}

\noindent
Here $(\mu_u)_{u\in\mathbb Z}$ follow the sine$_\beta$ point process, and $(g_u)_{u\in\mathbb Z}$ are i.i.d. standard Gaussians $\mathcal{N}_\beta(0,1)$ independent of $(\mu_u)$.

\end{theorem}

\bigskip
\begin{remark}[K-step consequence]
The limiting statistics of the $K$-step overlap matrix follow from the product (convolution) property
\[
\big\langle \mathbf{u}_i^{(N+K)},\, \mathbf{u}_j^{(N)} \big\rangle
= \Big( \prod_{s=1}^{K} \Omega^{(N+s)} \Big)_{i,j}
= \big( \Omega^{(N+K)} \cdots \Omega^{(N+1)} \big)_{i,j},
\]
where $\Omega^{(N)}$ denotes the one-step overlap matrix.
To characterize the $K$-step overlaps, we use the \emph{joint} law of the one-step overlaps. 
As in \cite{huang2022}, the asymptotic evolution of the eigenvectors is governed by a Markov chain; consequently, the $K$-step overlaps are obtained by multiplying the successive one-step overlap matrices.
Importantly, the matrices $\Omega^{(N+s)}$ are \emph{not} independent: each is an observable of the underlying bead process associated with the eigenvalue configuration, so their joint distribution is coupled through that process.
\end{remark}

\paragraph*{The eigenvector bead process.}
\noindent
We describe a Markov chain on a separable Hilbert space, indexed by a point process. Let $\mathcal{H}$ be a separable Hilbert space. We define the coupled Markov chain via the recurrence relations
\begin{equation}\label{eq:psi-phi-recursions}
\begin{aligned}
\Theta^{(s+1)} &= \Psi\!\left(\Theta^{(s)},\, G^{(s)},\, h\right),\quad
\Upsilon^{(s+1)} &= \Phi\!\left(\Upsilon^{(s)},\, \Theta^{(s+1)},\, \Theta^{(s)},\, G^{(s)}\right).
\end{aligned}
\end{equation}

where:
\begin{itemize}
\item $G^{(s)}=\{g_j^{(s)}\}_{j\in\mathbb{Z}}$, with $g_j^{(s)}\stackrel{\text{i.i.d.}}{\sim}\mathcal{N}_\beta(0,1)$, where the random variables are independent both across the index $j$ and across different values of $s$;
$h\in\mathbb{R}$; 
$\Theta$ denotes the space of locally finite point configurations on $\mathbb{R}$; 
$\{\Upsilon_j^{(s)}\}_{j\in\mathbb{Z}}$ is an orthogonal basis of $\mathcal{H}$.

\item \( \Psi \) is defined as:
\[
\Psi\!\left( \{ \Theta^{(s)}, G^{(s)}, h \right) 
= \left\{ z \in \mathbb{R} \, : \, \sum_{j \in \mathbb{Z}} \frac{|g_j^{(s)}|^2}{\mu_j - z} = h \right\}.
\]
\item \( \Phi \) is defined as:
\begin{align*}
\Phi\!\left(\Upsilon^{(s)},\, \Theta^{(s+1)},\, \Theta^{(s)},\, G^{(s)}\right)
&= \upsilon_i^{(s)} 
\sum_{j\in\mathbb{Z}} 
\frac{g_j^{(s)}}{\Theta_i^{(s+1)}-\Theta_j^{(s)}} 
\cdot \Upsilon_j^{(s)} \\[6pt]
\text{with} \quad
\upsilon_i^{(s)}
&= \left[
\sum_{j\in\mathbb{Z}}
\left(
\frac{|g_j^{(s)}|}{\Theta_i^{(s+1)}-\Theta_j^{(s)}}
\right)^2
\right]^{-1/2}
\end{align*}

\end{itemize}

\noindent
Fix an integer \( K \geq 1 \), and let \( E \in (-2, 2) \) be a given energy level.  
Let \( H_{N+K} \) be a Wigner matrix of size \( (N+K)\times(N+K) \). For each \( 0 \leq s \leq K \), let \( H_{N+K}^{(s)} \) denote the top-left \( (N+s)\times(N+s) \) principal minor of \( H_{N+K} \). Let \( \lambda_1^{(s)} \leq \cdots \leq \lambda_{N+s}^{(s)} \) be the ordered eigenvalues, and \( \mathbf{u}_1^{(s)}, \ldots, \mathbf{u}_{N+s}^{(s)} \) be the corresponding eigenvectors. Define the point process
\[
\Theta_N^{(s)} := \sum_{j=1}^{N+s} \delta_{\sqrt{4-E^2}\sqrt{N}\,(\lambda_j^{(s)} - E\sqrt{N})}.
\]
\noindent
For each \(s\ge 0\), define the index shift
\(
  i_0^{(s)} \coloneqq \min\bigl\{\, n \in \mathbb{N} : \lambda_n^{(s)} \ge E\sqrt{N} \,\bigr\}.
\) \\
Initialize
\(
\Upsilon_N^{(s=0)} \coloneqq \bigl(\mathbf{u}_{\,j+i_0^{(s=0)}}^{(s=0)}\bigr)_{j}
\) and define the entire sequence \( (\Upsilon_N^{(s)})_{s\ge 0} \) in accordance with
Conventions~\ref{conv:canonical-embedding}--\ref{conv:eigvec-normalization}
(see Section~\ref{sec:derivation-of-results}).

\smallskip
\noindent
\textbf{Conjecture.} 
The sequence
\(
\left( \Theta_N^{(s)}, \Upsilon_N^{(s)} \right)_{0 \leq s \leq K}
\)
defines a coupled Markov chain, which converges in distribution as \( N \to \infty \) 
to the Markov chain \eqref{eq:psi-phi-recursions} associated with the level
\[
h = -\frac{E}{2\sqrt{4-E^2}}.
\]

\subsection{Wishart Ensemble}

This section extends the analysis of eigenvector overlaps from the case of Wigner matrices to the setting of Wishart matrices, that is, empirical covariance matrices. It comes with notable structural differences: the eigenvalues of Wishart matrices follow the Marchenko--Pastur law \cite{marchenko1967}, instead of the semicircle law characteristic of Wigner matrices, and their construction as matrix products $XX^T$ introduces correlations between entries that are absent in the Wigner case.

Despite these challenges, recent works indicate that eigenvector overlaps exhibit a universal structure within the Wishart framework. Attal and Allez studied the overlaps between a large random covariance matrix and its principal submatrices, and derived explicit asymptotic expressions for the second-order moments of these overlaps as a function of the spectral position \cite{attal-allez2025}. Riabov, Tikhonov, and Bouchaud, on the other hand, investigated the case of sample covariance matrices constructed from a general population covariance matrix \(C\), focusing on the computation of the bi-resolvent of the model within this more general setting \cite{riabov2025}. These developments extend the theoretical framework established by earlier works on eigenvector dynamics under free addition \cite{allez-bouchaud2014} and on the structure of correlated random matrices \cite{bun-bouchaud-potters2018}.

We show that the structure of the limiting process — of the bead-type observed previously in the Wigner case — persists in the Wishart setting.

\paragraph{Wishart Matrices.}

Let $X = [x_{ij}]_{1 \le i \le N,\ 1 \le j \le T}$ be a real matrix whose entries are independent, centered, and standardized random variables, satisfying
\[
\mathbb{E}[x_{ij}] = 0, \qquad \mathbb{E}[x_{ij}^2] = 1, \qquad \text{and} \qquad \sup_{i,j}\, \mathbb{E}\!\left[|x_{ij}|^4\right] < \infty.
\]
We consider the asymptotic regime where both $N$ and $T$ tend to infinity while their ratio remains fixed,
\(
q = \frac{N}{T} \le 1.
\)
The Wishart matrix is defined by
\[
W_N = X^\top X \in \mathbb{R}^{N \times N}.
\]

\noindent
Fix a  integer $K \ge 0$. Let $W_{N+K}$ be a Wishart matrix as defined above.  
For any integer $0 \le s \le K$, we denote by
\(
W^{(s)}_{N+K} = [w_{ij}]_{1 \le i,j \le N+s}
\)
the top-left $(N+s) \times (N+s)$ principal minor of $W_{N+K}$,
with eigenvalues
\(
\lambda^{(s)}_1 \le \lambda^{(s)}_2 \le \cdots \le \lambda^{(s)}_{N+s},
\)
and corresponding normalized eigenvectors
\(
\mathbf{u}^{(s)}_1, \mathbf{u}^{(s)}_2, \ldots, \mathbf{u}^{(s)}_{N+s}.
\)
Then $W^{(s)}_{N+K}$ admits the spectral decomposition
\[
W^{(s)}_{N+K} = U^{(s)} \Lambda^{(s)} (U^{(s)})^\top,
\qquad
\Lambda^{(s)} = \mathrm{diag}\{\lambda^{(s)}_1, \lambda^{(s)}_2, \ldots, \lambda^{(s)}_{N+s}\},
\quad
U^{(s)} = [\mathbf{u}^{(s)}_1, \mathbf{u}^{(s)}_2, \ldots, \mathbf{u}^{(s)}_{N+s}].
\]

\medskip
\noindent
The empirical eigenvalue distribution of the rescaled matrix $\tfrac{1}{T} X^\top X$ converges, as $N,T \to \infty$, to the \emph{Marchenko--Pastur law} with density:
\begin{equation}
\rho_{\mathrm{MP}}(x)
= \frac{\sqrt{(\lambda_+ - x)(x - \lambda_-)}}{2\pi q x}\,
\mathbf{1}_{[\lambda_-,\,\lambda_+]}(x),
\quad
\lambda_{\pm} = (1 \pm \sqrt{q})^2,
\label{eq:MP-law}
\end{equation}

\begin{remark}
The nonzero eigenvalues of the matrices $X^{\top}X$ and $XX^{\top}$ are identical. 
Moreover, the eigenvectors of $XX^{\top}$ and $X^{\top}X$ associated with the same nonzero eigenvalue 
are mapped to each other by the matrix $X$. 
Indeed, if $v$ is an eigenvector of $X^{\top}X$ corresponding to a nonzero eigenvalue $\lambda$, 
then $u = Xv$ is an eigenvector of $XX^{\top}$ associated with the same $\lambda$. 
In our analysis, we frequently switch between the representations $X^{\top}X$ and $XX^{\top}$, 
since they share the same nonzero spectrum and their eigenvectors are directly related through this correspondence. 
For convenience, we therefore restrict our study to the case $q \le 1$, without loss of generality.
\end{remark}

\subsubsection{Eigenvalues for the Minors of Wishart Matrices}

\noindent
As in the work of Huang~\cite{huang2022}, the cases of eigenvalues located in the bulk and at the edge require separate analysis. We first state the results for the soft edges (both left and right), and then treat the hard-edge case \(q=1 \), where the dynamics resemble those observed in the bulk.

\paragraph{Eigenvalue Dynamics at the Soft Edge.}
\begin{theorem}[Left Edge]\label{thm:left-edge}
Fix integers $K, \ell \ge 1$, and let $q = N/T \in (0,1)$ be kept fixed in the asymptotic regime $N,T \to \infty$.
Consider a Wishart matrix $W_{N+K}$, and let $\{\lambda_j^{(s)}\}_{j}$ denote its increasingly ordered eigenvalues at level $s$. Then the array of rescaled left-edge increments
\[
\left(
-\frac{\sqrt{q}}{1-\sqrt{q}}\bigl(\lambda_{1}^{(s+1)} - \lambda_{1}^{(s)}\bigr),\,
-\frac{\sqrt{q}}{1-\sqrt{q}}\bigl(\lambda_{2}^{(s+1)} - \lambda_{2}^{(s)}\bigr),\,
\dots,\,
-\frac{\sqrt{q}}{1-\sqrt{q}}\bigl(\lambda_{\ell}^{(s+1)} - \lambda_{\ell}^{(s)}\bigr)
\right)_{0 \le s \le K-1}
\]
converges in distribution, as $N,T \to \infty$, to a random array with independent
Gamma-distributed entries on $(0,\infty)$ having density
\( \frac{1}{\sqrt{2\pi}}\, x^{-1/2} e^{-x/2}.
\)
Moreover, this array is independent of the array of the $\ell$ smallest eigenvalues,
\(
\left(
\lambda_{1}^{(0)},\,
\lambda_{2}^{(0)},\,
\dots,\,
\lambda_{\ell}^{(0)}
\right).
\)
\end{theorem}

\begin{theorem}[Right edge]\label{thm:right-edge}
Fix integers $K,\ell\ge 1$ and let $q=N/T\in(0,1]$ be kept constant in the asymptotic regime $N,T\to\infty$.
Consider a Wishart matrix $W_{N+K}$ and let $\{\lambda_j^{(s)}\}_j$ denote the (increasingly ordered) eigenvalues at level $s$.
Then the array of edge increments
\[
\left(
\frac{\sqrt{q}}{1+\sqrt{q}}\bigl(\lambda_{N+s+1}^{(s+1)}-\lambda_{N+s}^{(s)}\bigr),\,
\frac{\sqrt{q}}{1+\sqrt{q}}\bigl(\lambda_{N+s}^{(s+1)}-\lambda_{N+s-1}^{(s)}\bigr),\,
\dots,\,
\frac{\sqrt{q}}{1+\sqrt{q}}\bigl(\lambda_{N+s+1-\ell}^{(s+1)}-\lambda_{N+s-\ell}^{(s)}\bigr)
\right)_{0\le s\le K-1}
\]
converges in distribution, as $N,T\to\infty$, to a random array with independent
Gamma-distributed entries on $(0,\infty)$ with density
$\frac{1}{\sqrt{2\pi}}\, x^{-1/2} e^{-x/2}$. Moreover, this array is independent of the array of eigenvalues as stated in Theorem~\ref{thm:left-edge}.
\end{theorem}

\paragraph{Eigenvalue Dynamics at the Hard Edge.}

\noindent
Let $\mathcal{L}$ denote the space of locally finite point configurations on $\mathbb{R}$,
and set $\mathcal{G} := \mathbb{R}^{\mathbb{N}}$.
Define the operator
\[
\mathcal{T}:\;\mathcal{L} \times \mathcal{G} \times \mathbb{R}_{+} \longrightarrow \mathcal{L},
\qquad
(\Theta,\{g_j\},\xi)\longmapsto \mathcal{T}(\Theta,\{g_j\},\chi).
\]

\begin{equation}\label{eq:T-operator}
\mathcal{T}\!\big(\{\mu_j\}_{j\in\mathbb{N}},\,\{g_j\}_{j\in\mathbb{N}};\,\chi\big)
:= 
\left\{\, z\in\mathbb{R}\;:\;
\sum_{j\in\mathbb{N}}\frac{g_j^{2}}{\mu_j - z}
- \frac{\chi}{z} = 0 \right\}.
\end{equation}

\noindent
Given $\Theta^{(0)} \in \mathcal{L}$, we generate $\Theta^{(s)}$ iteratively by
\[
\Theta^{(s+1)} 
= \mathcal{T}\!\big(\Theta^{(s)},\,\{g^{(s)}_j\}_{j\in\mathbb{N}},\,\chi^{(s)}\big),
\qquad s = 0, 1, 2, \dots.
\]

\medskip
\noindent
\noindent
Fix integers \(K \ge 1\) and \(\alpha > K\).
For each $s \ge 0$, let
\(
\chi^{(s)} \sim \chi_{\alpha - s},
\)where $\chi_{\alpha-s}$ denotes a chi-square random variable with $\alpha-s$ degrees of freedom,
\(
\{ g^{(s)}_j \}_{j \in \mathbb{N}} \text{ i.i.d. } \mathcal{N}(0,1),
\). We assume that all random variables
\(
\{ \chi^{(s)} \}_{s \ge 0}
\)
and
\(
\{ g^{(s)}_j \}_{j \in \mathbb{N},\, s \ge 0}
\)
are mutually independent.

\begin{theorem}[Hard edge, \(q=1\)]\label{thm:hard-edge-q1}
Fix \(K\in\mathbb{N}\) and an integer \(\alpha > K\).
Let \(T - N = \alpha\).
For each \(s\in\{0,\dots,K\}\), let \(\{\lambda^{(s)}_{j}\}_{j=1}^{N+s}\) denote the spectrum of the matrix \(X_{N+s}^\top X_{N+s}\),
and define the point process
\[
\Theta_{N}^{(s)} \;=\; \sum_{j=1}^{N+s} \delta_{\,4N\,\lambda_{j}^{(s)}}.
\]
Then the process \(\big(\Theta_{N}^{(s)}\big)_{s=0}^{K}\) forms a Markov chain in \(s\) and converges in distribution to the Markov chain driven by the operator \(\mathcal{T}\).
The limiting point process of eigenvalues
\(
\Theta^{(s)} \;=\; \lim_{N\to\infty} \Theta_{N}^{(s)}
\)
is the \emph{Bessel point process} of parameter~\(\alpha\), describing the universal statistics at the hard edge of the spectrum. \cite{forrester1993hard,tracy1994bessel,forrester2010log}
\end{theorem}

\paragraph{Eigenvalue Dynamics in the Bulk.}

\begin{theorem}[Bulk regime]\label{thm:bulk}
Let $K \ge 1$ and fix an energy level
\(
E \in \bigl((1-\sqrt{q})^2,\,(1+\sqrt{q})^2\bigr), q = \frac{N}{T} \in (0,1),
\)
kept constant as $N,T \to \infty$. We recall the notation for $\rho_{\mathrm{MP}}$ (see Equation~\ref{eq:MP-law}).
Define the point process
\[
\Theta_N^{(s)}
\;=\;
\sum_{j=1}^{N+s}
\delta_{2\pi\rho_{\mathrm{MP}}(E)\,(\lambda_j^{(s)} - E T)},
\qquad 0 \le s \le K.
\]
Then $\bigl(\Theta_N^{(s)}\bigr)_{s=0}^K$ forms a Markov chain, and as $N,T\to\infty$,
it converges in distribution to the Markov chain
\(
\Theta^{(s+1)} = \mathcal{D}\bigl(\Theta^{(s)},\,\Gamma^{(s)},\,h\bigr)
\)
described in \cite[Eq.\,(1.4)]{huang2022},
with level parameter
\[
h(E, q)
= -\,\frac{q\,\big(E + q - 1\big)}{2\,\sqrt{\big((1+\sqrt{q})^{2} - E\big)\,\big(E - (1-\sqrt{q})^{2}\big)}}\,.
\]

\end{theorem}

\subsubsection{Eigenvectors for the Minors of Wishart Matrices}

\noindent
The results for the marginal laws are similar in the case of Wishart matrices, 
up to different constants. This further illustrates the universality of these processes. We also address the case of the hard edge.

\smallskip
\noindent
In the following, we study the eigenvalue gaps and the eigenvector overlaps.
For a unit eigenvector $\mathbf{u}^{(N+1)}_i$ of $W_{N+1}$ associated with the eigenvalue 
$\lambda^{(N+1)}_i$, define the overlaps $\Omega^{(N)}_{i,j}$ and the eigenvalue gaps 
$\Delta^{(N)}_{i,j}$ by
\[
\Omega^{(N)}_{i,j}
:= \big\langle \mathbf{u}^{(N+1)}_i,\, \mathbf{u}^{(N)}_j \big\rangle,
\qquad
\Omega^{(N)}_{i,N+1}
:= \big(\mathbf{u}^{(N+1)}_i\big)_{N+1},
\]
\[
\Delta^{(N)}_{i,j}
:= \lambda^{(N+1)}_i - \lambda^{(N)}_j,
\qquad
1 \le i \le N+1,\; 1 \le j \le N.
\]

\paragraph{Eigenvector Dynamics at the Soft Edge.}

\noindent
The classical results of Johnstone~\cite{johnstone2001} and Soshnikov~\cite{soshnikov2002} establish the universal nature of the soft–edge behavior in Wishart matrices. Johnstone showed that, for sample covariance matrices, the largest eigenvalue—after centering at $\lambda_+ = (1+\sqrt{q})^2$ and rescaling by $N^{2/3}$—converges to the Tracy–Widom law. Soshnikov extended this result by proving that the same limiting distribution holds for a broad class of population distributions, thereby confirming the universality of the Airy process at the spectral edge. These works form the theoretical basis for the asymptotic analysis stated in Theorem~\ref{thm:soft-edge-case}.

\begin{theorem}[Soft Edge]\label{thm:soft-edge-case}
Fix an integer $\ell \ge 1$.
For $1 \le i,j \le \ell$, the overlap matrix converges to the identity matrix, that is,
\[
\lim_{N \to \infty} \bigl(\Omega^{(N)}_{i,j}\bigr)_{1 \le i,j \le \ell} = I_\ell.
\]
Moreover, there exists a constant $c_q > 0$, depending only on the aspect ratio $q$, such that the leading-order correction scales as $O(N^{-1/3})$ and satisfies
\begin{equation}\label{eq:soft-edge-edge}
  N^{1/3}\bigl(\Omega^{(N)}_{i,j} - (I_\ell)_{i,j}\bigr)
  \;\underset{N \to \infty}{\overset{d}{\sim}}\;
  \begin{cases}
    c_q\,\dfrac{g_i g_j}{\alpha_j - \alpha_i}, & i \ne j,\\[0.8ex]
    0, & i = j,
  \end{cases}
  \qquad 1 \le i,j \le \ell.
\end{equation}
Here \((\alpha_i)_{i=1}^\ell\) follow the Airy$_{\beta}$ point process, and \((g_i)_{i=1}^\ell\) are i.i.d.\ standard Gaussians \(\mathcal{N}(0,1)\) independent of \((\alpha_i)\). For the right edge, set \(c_q = (1+\sqrt{q})^{-1/3}\) and for the left edge, set \(c_q = (1-\sqrt{q})^{-1/3}\).
\end{theorem}

\paragraph{Eigenvector Dynamics in the Bulk.}
In the bulk of Wishart ensembles, the rescaled eigenvalues converge to the \textit{sine}$_\beta$ point process. This universal limit, first identified by Dyson and Mehta~\cite{dyson1962statistical} and detailed by Forrester~\cite{forrester2010log}, describes the local correlations of eigenvalues governed by the sine kernel ($\beta=2$) or its Pfaffian analogue ($\beta=1$).

\begin{theorem}[Bulk, real Wishart]\label{thm:bulk-case-loe}
We recall Convention~\ref{conv:inverse-stieltjes-branches}. 
Fix $E \in \bigl((1-\sqrt{q})^2,(1+\sqrt{q})^2\bigr)$ and define the anchor index
\[
i_N(E) \;:=\; \min\bigl\{\,k \in \{1,\dots,N\}\ :\ \lambda^{(N)}_{k} \ge E\,T\,\bigr\}.
\]
We recall the notation for the Marčenko–Pastur density $\rho_{\mathrm{MP}}$ (see Equation~\ref{eq:MP-law}).
Then the overlaps satisfy
\begin{equation}\label{eq:bulk-near-anchored-loe}
\Omega^{(N)}_{\,i_N(E)+u,\;i_N(E)+v}
\;\mathrel{\underset{N\to\infty}{\overset{d}{\sim}}}\;
\frac{g_v}{\,S_{\mu,\,g^2}^{-1}(h)_u - \mu_v\,}\,
\sqrt{\bigl(S_{\mu,\,g^2}^{-1}\bigr)'(h)_u},
\qquad 
h
\;=\;
-\,\frac{E + q - 1}{2E}
\cdot
\frac{1}{2\pi\rho_{\mathrm{MP}}(E)}\,.
\end{equation}

\noindent
Here $(\mu_u)_{u \in \mathbb{Z}}$ follows the sine\(_1\) point process, 
and $(g_u)_{u \in \mathbb{Z}}$ are i.i.d.\ real Gaussian random variables 
$g_u \sim \mathcal{N}(0,1)$, independent of $(\mu_u)$.
\end{theorem}

\noindent
\textbf{Conjecture.} 
The coupled eigenvalue--eigenvector dynamics of Wishart matrices 
are expected to converge to the same limiting Markov chain as in the Wigner case, 
with the parameter \( h \) specified in Theorem~\ref{thm:bulk}.

\paragraph{Eigenvector Dynamics at the Hard Edge.}
At the hard edge of Wishart ensembles, the smallest eigenvalues—after appropriate rescaling—converge to the Bessel point process. This limit was first derived by Tracy and Widom~\cite{tracy1994bessel}, who expressed the distribution of the smallest eigenvalues in terms of a Fredholm determinant with the Bessel kernel. Forrester~\cite{forrester1993hard,forrester2010log} further developed the analysis of the corresponding correlation functions and confirmed the universality of this behavior for general $\beta$–ensembles. These results form the foundation for the hard–edge asymptotics stated in Theorem~\ref{thm:hard-edge-case-wishart}.

\begin{theorem}[Hard Edge]\label{thm:hard-edge-case-wishart}
We recall Convention~\ref{conv:inverse-stieltjes-branches}, and use the same convention for the inverse branches of the Stieltjes transform appearing below. 
Let $\alpha = T - N \ge 0$ denote the rectangularity parameter of the Wishart ensemble.
For $u, v \in \mathbb{N}$, the overlaps satisfy
\begin{equation}\label{eq:hard-edge-overlaps}
\Omega^{(N)}_{u,v}
\;\mathrel{\underset{N \to \infty}{\overset{d}{\sim}}}\;
\frac{g_v \, \sqrt{\xi_v}}{\,D_{g^2, \xi, \alpha}^{-1}(0)_u - \xi_v\,}
\left(
    \sum_{v' \in \mathbb{N}} 
    \frac{(g_{v'})^2\xi_{v'}}{\bigl(D_{g^2, \xi, \alpha}^{-1}(0)_u - \xi_{v'}\bigr)^2}
\right)^{-1/2}.
\end{equation}

\noindent
Here $(\xi_u)_{u \in \mathbb{N}}$ follow the Bessel point process associated with the parameter $\alpha$, 
and $(g_u)_{u \in \mathbb{N}}$ are i.i.d.~standard Gaussian random variables 
$\mathcal{N}(0,1)$, independent of $(\xi_u)$.

\medskip
\noindent
The function $D_{g^2, \xi, \alpha}(z)$ is defined by
\begin{equation}\label{eq:def-D-xi}
D_{g^2, \xi, \alpha}(z)
= \sum_{i \in \mathbb{N}} \frac{g_i^2}{ \xi_i-z}
- \frac{\chi_{\alpha}}{z},
\end{equation}
where $\chi_{\alpha}$ denotes a chi-square random variable with $\alpha$ degrees of freedom.
\end{theorem}

\subsection{Numerical Simulations}

\subsubsection{Wigner Ensemble}
\begin{figure}[H]
  \centering

  \begin{subfigure}{0.45\linewidth}
    \includegraphics[width=0.8\linewidth]{\detokenize{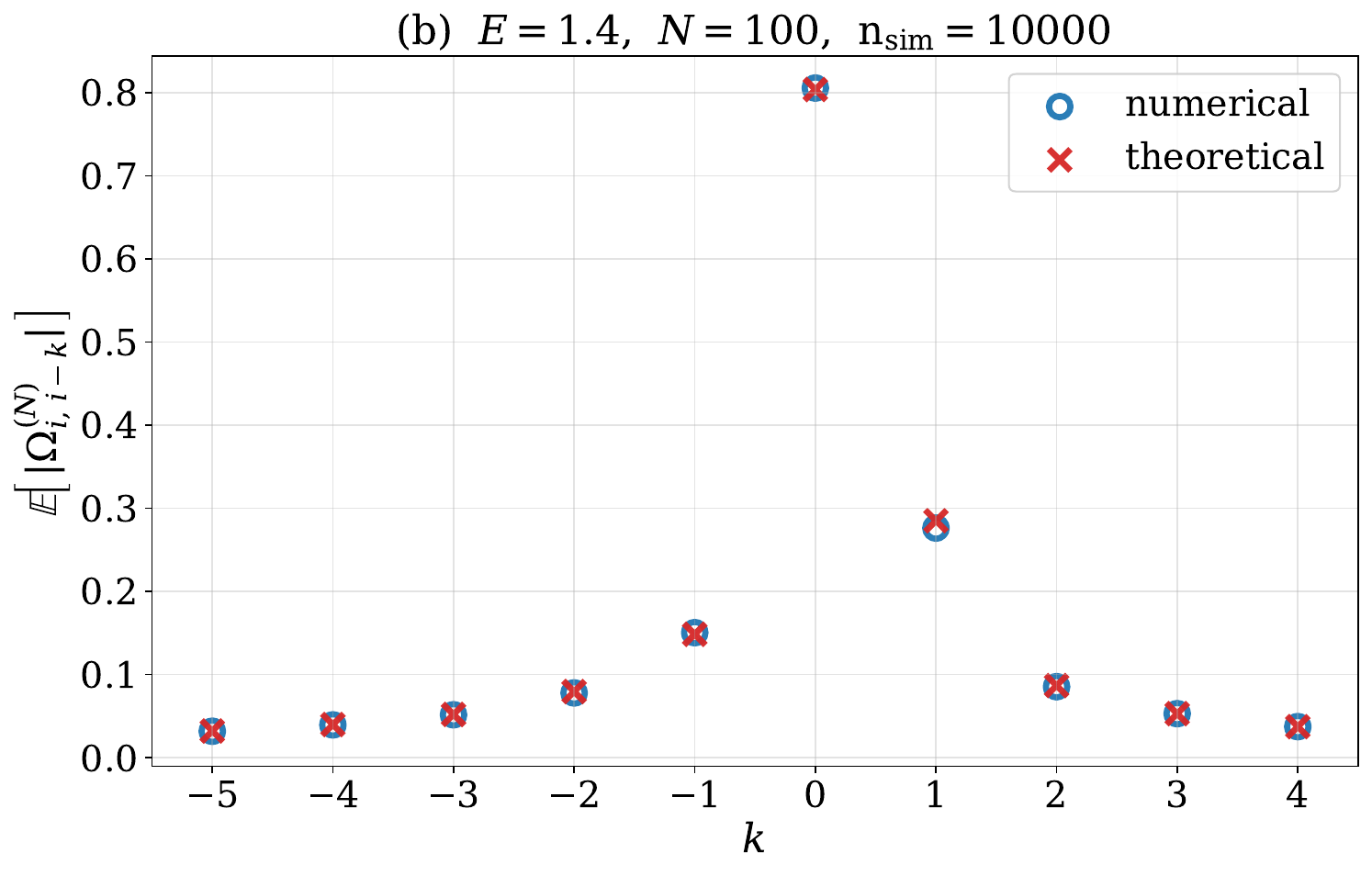}}
  \end{subfigure}\hfill
  \begin{subfigure}{0.45\linewidth}
    \includegraphics[width=0.8\linewidth]{\detokenize{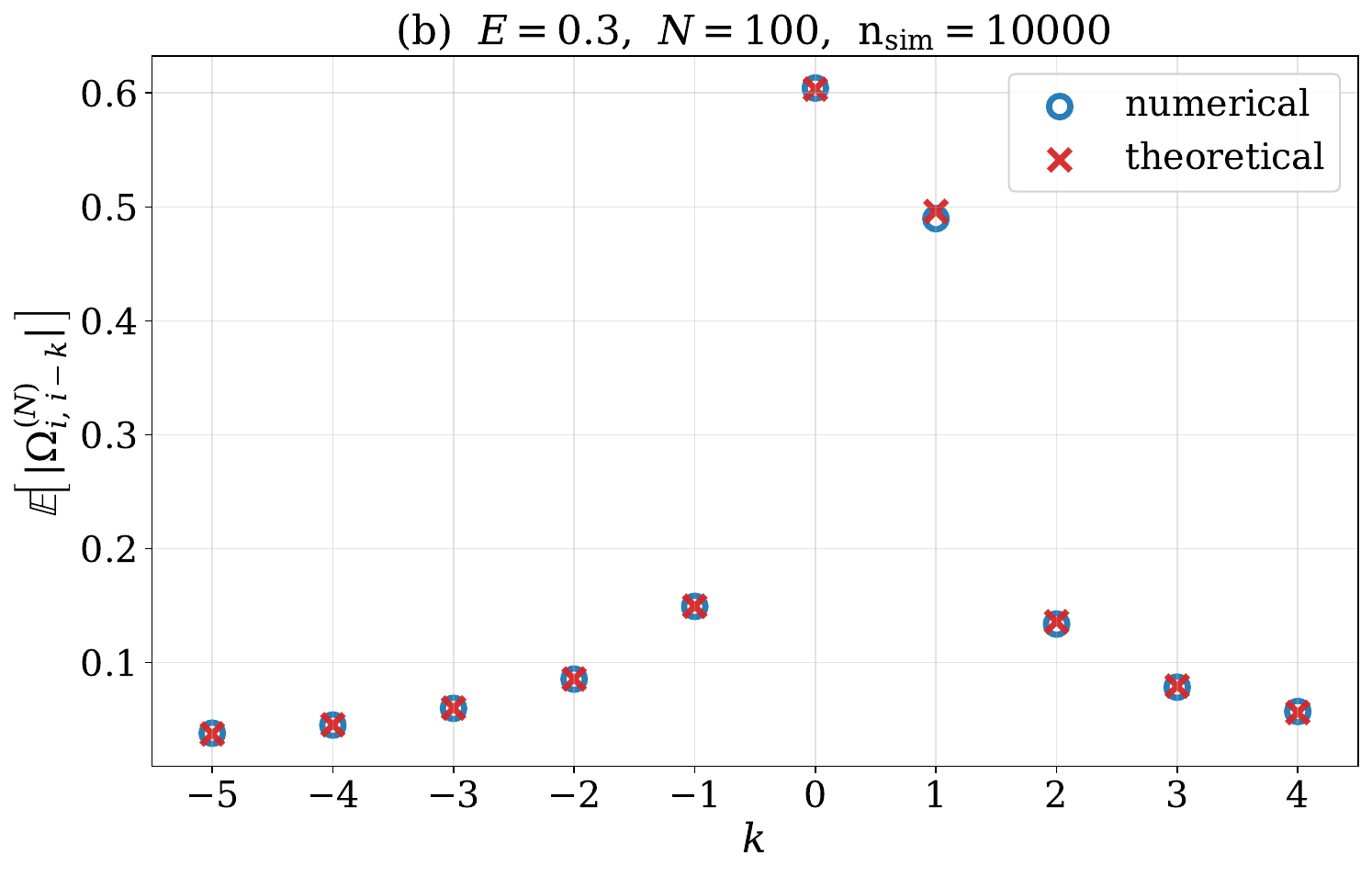}}
  \end{subfigure}

\caption{Mean absolute overlap as a function of the offset $k=i-j$. Panels (a)–(b) plot $\mathbb{E}\left[|\Omega^{(N)}_{i,i-k}|\right]$ for $N=100$ at energies (a) $E=1.4$ and (b) $E=0.3$. Blue circles show the \emph{empirical} averages: in each run we draw a \emph{real-symmetric Wigner} matrix ($\beta=1$), augment it by one row and one column, and then compute the overlap between the relevant eigenvectors. Red dashed squares show the \emph{theoretical} simulation: we sample the \emph{sine--$\beta$} process at $\beta=1$ with Gaussian variables, and evolve the inserted level via the secular equation $S(z)=h$. Evaluating the derivative $S'(z)$ at the corresponding root yields the predicted overlaps according to formula~\eqref{eq:bulk-near-anchored}. Results are averaged over $10{,}000$ Monte Carlo runs.}

  \label{fig:overlap-means}
\end{figure}

\begin{figure}[H]
  \centering
  \includegraphics[width=0.7\linewidth]{\detokenize{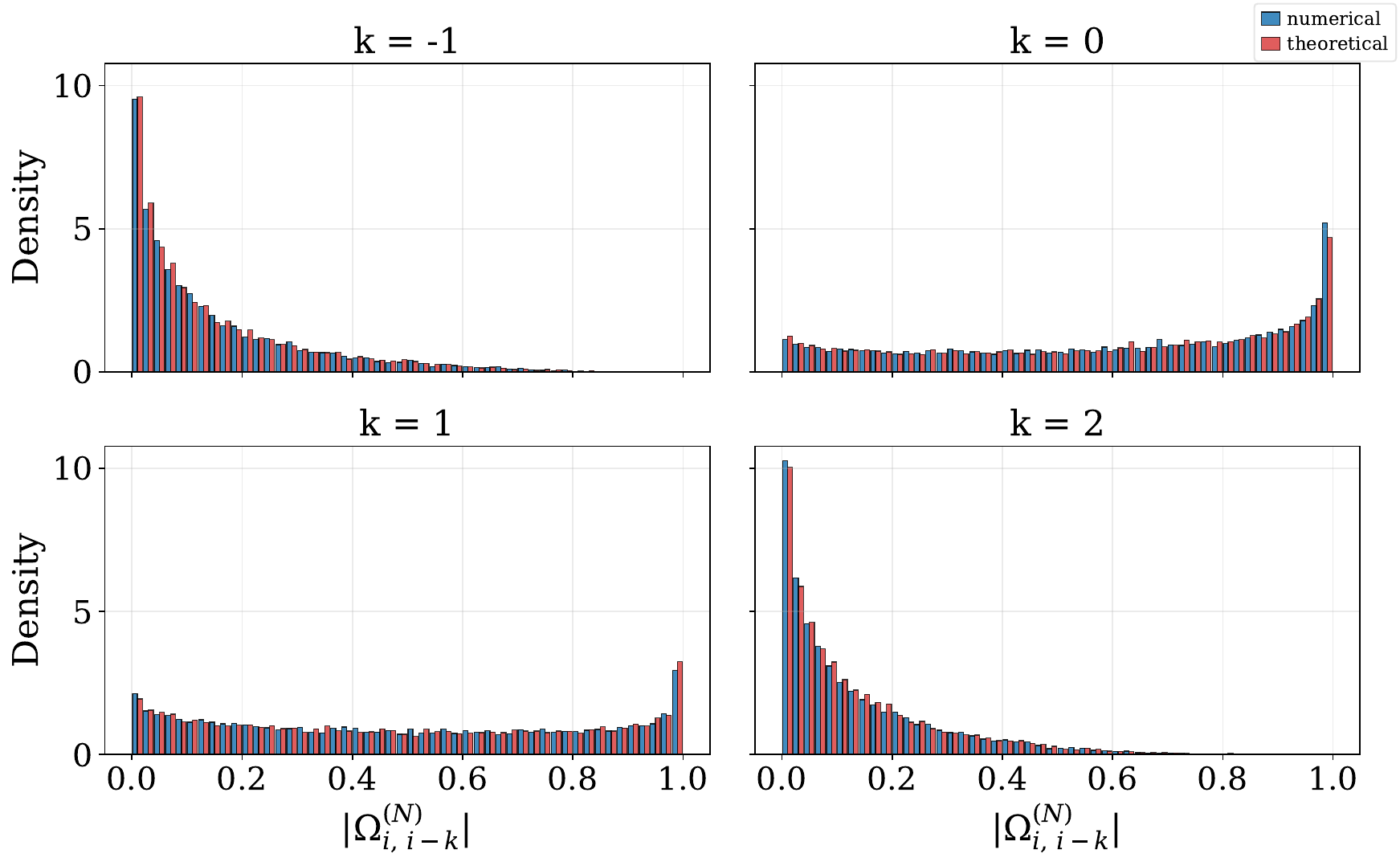}}
\caption{Distribution of absolute overlaps at offsets $k\in\{-1,0,1,2\}$. Panels compare the empirical (\textit{numerical}, blue) and predicted (\textit{theoretical}, red) densities of $|\Omega^{(N)}_{i,\,i-k}|$ for $N=100$ at energy $E=0.3$. We fix the reference index $i=60$ and examine the four neighbour indices $j=i+k\in\{58,59,60,61\}$ (i.e., $k=-1,0,1,2$). Each panel aggregates $10{,}000$ Monte Carlo runs. For the empirical distributions, in each run we draw a \emph{real-symmetric Wigner} matrix ($\beta=1$), append one row and one column to form an $(N{+}1)\times(N{+}1)$ arrowhead matrix, and read the overlaps from its eigenvectors. For the theoretical distributions, we sample the \emph{sine--$\beta$} process at $\beta=1$ with Gaussian variables, and evolve the inserted level via the secular equation $S(z)=h$. Evaluating the derivative $S'(z)$ at the corresponding root yields the predicted overlaps according to formula~\eqref{eq:bulk-near-anchored}. Histograms share a common binning and are normalized to unit area.
}
\label{fig:overlaps-hist-k--2--1-0-1}

\end{figure}

\noindent
As shown in Figure~1, the projection of the new eigenvector onto the previous eigenvector at the same energy level exhibits a sharply peaked shape. The width of this peak is independent of \(N\) and tends to become increasingly sharp as \(|E|\) grows. This is consistent with the description of \(h = -\frac{E}{2\sqrt{4-E^2}}\) given in \cite{huang2022}.
\noindent
In Figure~2, for \(k=1\) one observes a mass at \(1\), corresponding to realizations in which the overlap jumps directly to the next level; this occurs when the Gaussian ratio \(g_{i+1}/g_i\) (which is Cauchy–distributed) takes a large value. Overall, the overlap distribution for \(k=1\) is shifted toward larger values compared with \(k=-1\).

\subsubsection{Wishart Ensemble}

\begin{figure}[H]
  \centering

 \includegraphics[width=0.9\linewidth]{\detokenize{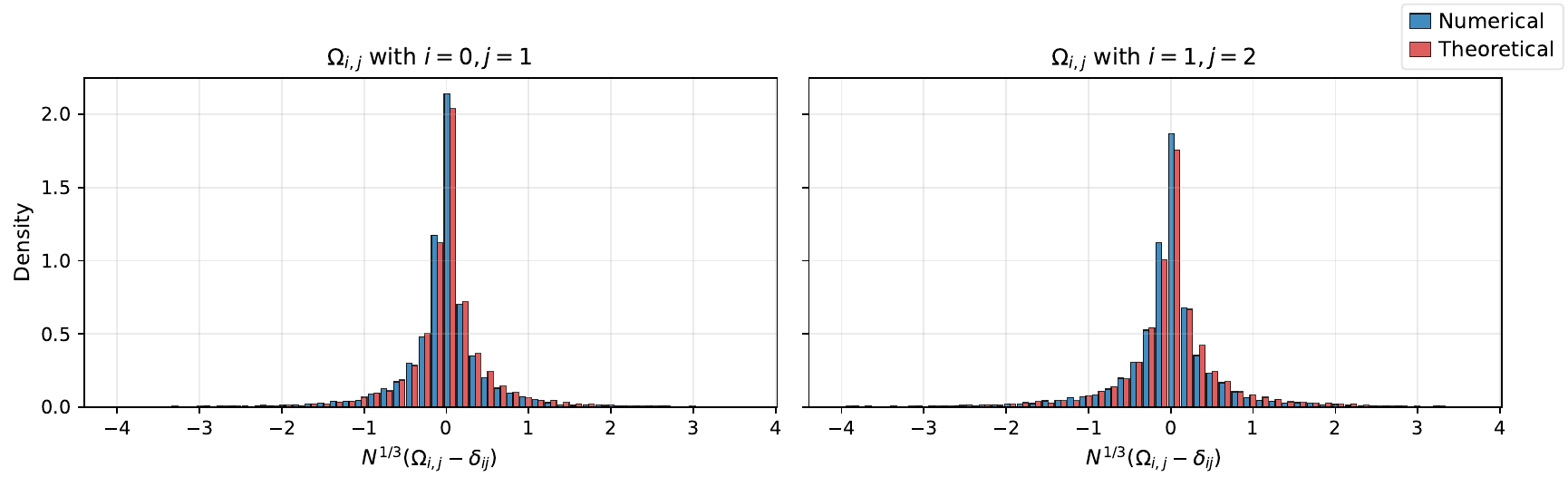}}

\caption{Distributions of rescaled overlaps at the right edge for $q=0.9$. Panels show the probability density of $N^{1/3}(\Omega_{i,j}-\delta_{ij})$ for $N=100$ with (left) $(i,j)=(0,1)$ and (right) $(i,j)=(1,2)$. Blue bars show the \emph{numerical} distributions: in each run we draw a \emph{real Wishart} matrix with aspect ratio $q=T/N$, augment it by one row and one column, and compute the overlap between the perturbed and unperturbed eigenvectors near the largest eigenvalue. Red bars show the \emph{theoretical} predictions: we sample the \emph{Airy--$\beta$} point process at $\beta=1$ (approximated via GOE), draw independent Gaussian variables $(g_i,g_j)$, and compute the overlaps according to the soft-edge formula in Theorem~\ref{thm:soft-edge-case} with constant $c_q=(1+\sqrt{q})^{-1/3}$. The agreement validates the soft-edge universality for Wishart matrices. Results are averaged over $10{,}000$ Monte Carlo simulations.}

  \label{fig:soft-case-distributions}
\end{figure}

\begin{figure}[H]
  \centering

 \includegraphics[width=0.8\linewidth]{\detokenize{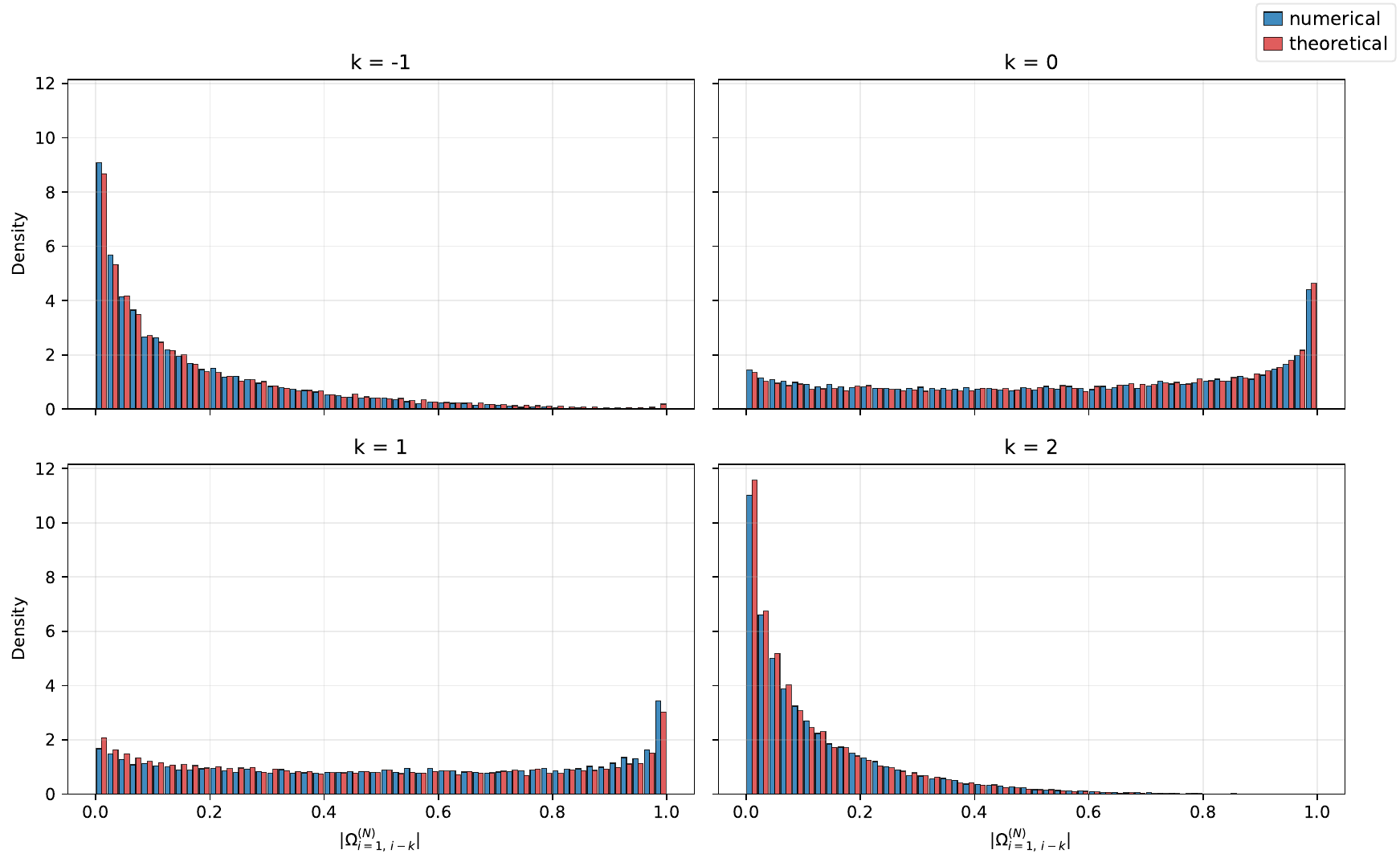}}

\caption{Distributions of rescaled overlaps at the hard edge for $\alpha=1$. Panels show the probability density of $|\Omega^{(N)}_{i,i-k}|$ for $N=200$ with offset values (top left) $k=-1$, (top right) $k=0$, (bottom left) $k=1$, and (bottom right) $k=2$, corresponding to overlaps between the smallest eigenvectors. Blue bars show the \emph{numerical} distributions: in each run we draw a \emph{real Wishart} matrix with rectangularity parameter $\alpha=T-N=1$, augment it by one row and one column, and compute the overlap between the perturbed and unperturbed eigenvectors near the smallest eigenvalue. Red bars show the \emph{theoretical} predictions: we sample the \emph{Bessel point process} at $\beta=1$ (approximated via small Wigner eigenvalues with $N_w=1000$, keeping the smallest $80$ points), draw independent Gaussian variables $(g_i)$, and compute the overlaps according to the hard-edge formula in Theorem~\ref{thm:hard-edge-case-wishart} by solving the secular equation $D_{g^2,\xi,\alpha}(z)=0$ where $D_{g^2,\xi,\alpha}(z)=\sum_i g_i^2/(\xi_i-z)-\chi^2_\alpha/z$. Results are averaged over $15{,}000$ Monte Carlo simulations.}

  \label{fig:hard-case-distributions}
\end{figure}

\noindent
The overlap distributions at the soft and hard edges of Wishart matrices shown in Figures~\ref{fig:soft-case-distributions} and~\ref{fig:hard-case-distributions} exhibit qualitatively similar behavior to their Wigner counterparts, confirming edge universality across ensembles.

\section{Derivation of the results}
\label{sec:derivation-of-results}

\medskip
\noindent
\textbf{Gaussian distributions.}
For $\beta\in\{1,2\}$ we write $\mathcal N_{\beta}(0,1)$ for the \emph{standard Gaussian of symmetry $\beta$}:
\[
\mathcal N_{1}(0,1)=\text{the real } \mathcal N_{\mathbb R}(0,1),\qquad
\mathcal N_{2}(0,1)=\text{the complex } \mathcal N_{\mathbb C}(0,1),
\]
where a complex standard Gaussian is defined by
\[
Z=X+iY,\qquad X,Y\stackrel{\text{i.i.d.}}{\sim}\mathcal N_{\mathbb R}\!\Big(0,\tfrac12\Big),\quad
\mathbb E[Z]=0,\quad \mathbb E[|Z|^2]=1 .
\]

\noindent The limiting eigenvalue distribution of Wigner matrices
(see Paragraph~\ref{par:wigner-matrices})
is given by the semicircle law, \[ \rho_{\mathrm{sc}}(x)=\tfrac{\sqrt{4-x^{2}}}{2\pi}, \; x\in[-2,2], \qquad \mu_{\mathrm{sc}}[a,b]=\int_a^b \rho_{\mathrm{sc}}(x)\,dx, \; [a,b]\subset[-2,2]. \]

\medskip
\noindent
\textbf{Asymptotic notation.}
We write $a_N=O(b_N)$ if $|a_N|\le C\,b_N$ for some $C>0$ (uniformly in $N$),
and $a_N=o(b_N)$ if $a_N/b_N\to 0$. ``With high probability'' (w.h.p.) means probability $1-o(1)$.
We write $a_N=O_{\mathbb P}(b_N)$ (resp.\ $o_{\mathbb P}(b_N)$) if the bound holds with high probability.

\medskip
\begin{convention}\label{conv:canonical-embedding}
For each \(N\ge1\), we identify \(\mathbb{R}^{N}\) with the coordinate subspace
\(\mathbb{R}^{N}\times\{0\}\subset\mathbb{R}^{N+1}\) via the canonical isometric embedding
\[
\iota_{N}:\mathbb{R}^{N}\hookrightarrow\mathbb{R}^{N+1},\qquad \iota_{N}(x)=(x,0).
\]
When no confusion arises we still write \(x\) for \(\iota_{N}(x)\).
In particular, the Euclidean inner product is preserved and \(e_{N+1}\) is orthogonal to \(\iota_{N}(\mathbb{R}^{N})\).
As a consequence, the unit sphere \(S^{N-1}\subset\mathbb{R}^{N}\) is canonically included in \(S^{N}\subset\mathbb{R}^{N+1}\) via \(x\mapsto(x,0)\).
\end{convention}

\medskip
\begin{convention}[Eigenvector sign/phase normalization]\label{conv:eigvec-normalization}
In the real case, we fix the sign ambiguity by imposing
\[
\Omega^{(N)}_{i,i} \ge 0 \quad \text{for all } i=1,\dots,N+1.
\]
In the complex case, eigenvectors are defined up to a global phase \(e^{\mathrm i\theta}\).
We fix this ambiguity by requiring that a distinguished coordinate be real and nonnegative; concretely,
\[
\Omega^{(N)}_{i,i}\in[0,\infty)\subset\mathbb{R}
\quad\text{for all } i=1,\dots,N+1,
\]
i.e., we multiply each eigenvector by a phase so that its \((N+1)\)-st entry is real and \(\ge 0\).
If this entry vanishes, choose the smallest index \(m\) with \((\Omega_i^{(s)})_m\neq 0\) and impose \((\Omega^{(N)}_{i})_m\in[0,\infty)\).
Equivalently,
\[
\arg\!\big\langle e_{N+1},\,\Omega^{(N)}_{i}\big\rangle=0
\quad\text{and}\quad
\big\langle e_{N+1},\,\Omega^{(N)}_{i}\big\rangle\ge 0.
\]
\end{convention}

\noindent
We now turn to the derivation of the results presented above. Our approach begins with a finite-\( N \) analysis, followed by an asymptotic investigation in both the spectral edge and bulk regimes. The standard strategy for studying the eigenvalues of macroscopic minors of Wigner matrices consists in analyzing the determinant equation at finite \( N \), and then passing to the limit. We shall adopt the same methodology.

\subsection{Master Equations for Eigenvalues and Eigenvectors in the Wigner case}

\begin{proposition}[Secular equation for the eigenvalues]\label{prop:secular}
If $g_j=0$ for some $j$, then $\lambda^{(N)}_j$ remains an eigenvalue of $H_{N+1}$ with eigenvector
$(u^{(N)}_j,0)$. All other eigenvalues of $H_{N+1}$ are the simple zeros of the secular function
\begin{equation}\label{eq:secular}
f(z):=g_{N+1}-z-\sum_{j=1}^N \frac{|g_j^{(N)}|^2}{\lambda^{(N)}_j-z}
\end{equation}
\end{proposition}

\begin{proposition}[Master equation for the eigenvectors]\label{prop:master-vec}
Assume the spectrum $\{\lambda^{(N)}_j\}_{j=1}^N$ is simple.
Let $\lambda^{(N+1)}_i$ be any zero of \eqref{eq:secular}. Then the overlaps satisfy
\begin{subequations}\label{eq:master-system}
\begin{align}
\Omega^{(N)}_{i,j}
&= \frac{1}{\sqrt{\bigl|f'(\lambda^{(N+1)}_i)\bigr|}}\,
   \frac{g_j^{(N)}}{\Delta_{i,j}^{(N)}},
\qquad j=1,\dots,N, \label{eq:master}\\[4pt]
f'(\lambda)
&= -\left( 1 + \sum_{j=1}^N
           \left(\frac{|g_j^{(N)}|}{\lambda-\lambda^{(N)}_j}\right)^{\!2} \right).
\label{eq:fprime}
\end{align}
\end{subequations}
\end{proposition}

\begin{remark}[Non-simple spectrum]
If $H_N$ has multiple eigenvalues, choose common eigenpairs for $H_N$ and $H_{N+1}$ within each shared invariant eigenspace.
On those subspaces the overlap matrix is the identity; the remaining overlaps obey \eqref{eq:master}.
\end{remark}

\begin{lemma}[Ratio identities and diagonal normalization]\label{lem:ratios}
For $g_i,g_j\neq0$ and any $k$,
\begin{equation}\label{eq:ratios}
\frac{\Omega^{(N)}_{k,j}}{\Omega^{(N)}_{k,i}}
=\frac{\Delta_{k,i}^{(N)}}{\Delta_{k,j}^{(N)}}\cdot \frac{g_j^{(N)}}{g_i^{(N)}},
\qquad
\frac{\Omega^{(N)}_{k,N+1}}{\Omega^{(N)}_{k,i}}
=\frac{\Delta_{k,i}^{(N)}}{g_i^{(N)}}.
\end{equation}
In addition, for $k=i$,
\begin{equation}\label{eq:diag-overlap}
\bigl(\Omega^{(N)}_{i,i}\bigr)^2
=
\left[
\left(\frac{\Delta_{i,i}^{(N)}}{|g_i^{(N)}|}\right)^{\!2}
\left(1+\sum_{j=1}^N \left(\frac{|g_j^{(N)}|}{\Delta_{i,j}^{(N)}}\right)^{\!2}\right)
\right]^{-1}.
\end{equation}
\end{lemma}

\begin{remark}[Cauchy law for overlap ratios]
Equation~\eqref{eq:ratios} exhibits the random ratio \(g_j/g_i\) with \(g_i,g_j\stackrel{\mathrm{i.i.d.}}{\sim}\mathcal N_\beta(0,1)\).
\noindent
Since the ratio of two independent centered Gaussian variables with the same variance is standard Cauchy, we have
\(
g_j/g_i\ \sim\ \mathrm{Cauchy}(0,1).
\)
The second factor \(\Delta_{k,i}/\Delta_{k,j}\) is a ratio of eigenvalue gaps. For distinct \(i\neq j\), the Cauchy interlacing theorem for rank-one extension constrains these gaps (e.g., their signs and sizes relative to the local spacing), so \(\Omega^{(N)}_{k,j}/\Omega^{(N)}_{k,i}\) is the product of a standard Cauchy variable and a factor driven by interlacing. The two factors are not independent.
\end{remark}

\begin{proposition}[Shift identity for symmetric perturbations]\label{prop:shift}
Let \(H,\mathcal{D}\in\mathcal S_N(\mathbb R)\), with normalized eigenpairs
\(\{(\lambda^H_j,\mathbf{u}^H_j)\}\) of \(H\) and \(\{(\lambda^{H+\mathcal{D}}_i,\mathbf{u}^{H+\mathcal{D}}_i)\}\) of \(H+\mathcal{D}\).
Whenever the denominator is nonzero,
\begin{equation}\label{eq:shift}
\lambda^{H+\mathcal{D}}_i-\lambda^H_j
=\frac{(\mathbf{u}^{H+\mathcal{D}}_i)^\top \mathcal{D}\,\mathbf{u}^H_j}{(\mathbf{u}^{H+\mathcal{D}}_i)^\top \mathbf{u}^H_j}.
\end{equation}
\end{proposition}

\noindent\emph{Outline.} The statements then follow in a cascade:
\eqref{eq:shift} $\Rightarrow$ \eqref{eq:ratios} (by forming a ratio to cancel the common factor)
$\Rightarrow$ \eqref{eq:diag-overlap} (by normalization)
$\Rightarrow$ \eqref{eq:master} (by differentiating the secular equation~\eqref{eq:secular}).

\begin{proof}[Proof]
We first derive the secular equation \eqref{eq:secular} via the Schur complement of the block matrix, which yields a polynomial equation in $\lambda$. Combining the eigen-equations $(H+\mathcal{D})\,\mathbf{u}_i^{H+\mathcal{D}}=\lambda_i^{H+\mathcal{D}}\mathbf{u}_i^{H+\mathcal{D}}$ and $H\,\mathbf{u}_j^H=\lambda_j^H\mathbf{u}_j^H$, taking the inner product with $\mathbf{u}_j^H$, and subtracting gives the shift identity
\(
(\lambda_i^{H+\mathcal{D}}-\lambda_j^H)\,\langle \mathbf{u}_j^H,\mathbf{u}_i^{H+\mathcal{D}}\rangle
=\langle \mathbf{u}_j^H,\mathcal{D}\,\mathbf{u}_i^{H+\mathcal{D}}\rangle,
\)
valid for all $i,j$. 
Applying \eqref{eq:shift} with $\mathcal{D}=\mathbf{e}_{N+1}\mathbf{g}^\top+\mathbf{g}\mathbf{e}_{N+1}^\top$ (supported on the last row/column) and then forming a ratio to cancel the common factor $(\mathbf{u}_i^{(N+1)})^\top \mathbf{e}_{N+1}$ yields~\eqref{eq:ratios}.
The normalization \(\sum_{j=1}^{N+1}(\Omega^{(N)}_{i,j})^2=1\) gives \eqref{eq:diag-overlap}.
Finally, \eqref{eq:master} is obtained by differentiating \eqref{eq:secular} at a zero \(\lambda\) and comparing with \eqref{eq:diag-overlap}.
\end{proof}

\subsection{Asymptotic Behavior of the Master Equation in the Wigner case}

\noindent
We investigate the asymptotic behavior of these relations both in the bulk and at the edge of the spectrum.

\subsubsection{Edge regime}
\noindent
In this section, we establish the behavior of the overlap matrix in regimes involving the spectral edges. We begin by recalling three ingredients from \cite{huang2022} and \cite{bao2025} that will be used to analyze the asymptotic behavior at the edge. 

\noindent

\begin{convention}[Working at the left edge]\label{conv:left-soft-edge}
We work throughout at the \emph{left} (lower) edge \(-2\sqrt{N}\).
The \emph{right} edge case is entirely analogous: replace
\(i \mapsto N-i\) and \(-2\sqrt{N} \mapsto +2\sqrt{N}\)
in all statements below.
\end{convention}

\begin{itemize}
\item \emph{Gaussian edge coordinates and asymptotic independence} \cite{huang2022}.
For any fixed $\ell\ge 1$ and as $N\to\infty$,
\begin{equation}\label{eq:gaussian-edge-coords}
  (g_1^{(N)},\dots,g_\ell^{(N)})\ \;\xrightarrow{d}\; \mathcal N_{\beta}(0,I_\ell).
\end{equation}
Equivalently, the coordinates $\{g_i^{(s)}\}_{i=1}^{\ell}$ are asymptotically i.i.d.\ standard Gaussians and are asymptotically independent of the (rescaled) edge eigenvalues.

  \item \emph{Edge gap \cite{huang2022}.}
  Let \(\lambda_i^{(N+1)}\) and \(\lambda_i^{(N)}\) be the \(i\)-th smallest eigenvalues of \(H_{N+1}\) and \(H_{N}\). Then
  \begin{equation}\label{eq:soft-edge-gap}
    \sqrt{N}\,\bigl(\lambda_i^{(N)}-\lambda_i^{(N+1)}\bigr)
    \;=\; \lvert g_i^{(N)}\rvert^{2} + o_{\mathbb P}(1).
  \end{equation}
In particular, \(|g_i^{(N)}|^{2}\) converges in distribution to a Gamma law with unit mean and variance \(2/\beta\),
independently of \(\lambda_i^{(N)}\); hence the scaled gap inherits this distributional limit.

\item \emph{Adjacent-step eigenvector alignment \cite{bao2025}.}
The central overlap is asymptotically aligned and
\begin{equation}\label{eq:adjacent-alignment}
  \big(\mathbf{u}_i^{(N+1)}\big)^\top \mathbf{u}_i^{(N)}
  \;=\; 1 \;+\; O_{\mathbb{P}}\!\left(N^{-2/3}\right).
\end{equation}

\end{itemize}

\paragraph*{Proof of Theorem~\ref{thm:edge-case}}\mbox{}\\

\noindent
Fix $\ell\ge 1$. Let $\{(\lambda_i^{(N)},\mathbf{u}_i^{(N)})\}_{i=1}^{N}$ be the eigenpairs of $H_{N}$ in increasing order, and set
\[
U_\ell^{(N)} := \big[\,\mathbf{u}_1^{(N)}\ \cdots\ \mathbf{u}_\ell^{(N)}\,\big]\in\mathbb R^{N\times \ell},
\]
the matrix of the $\ell$ \emph{smallest} (left-edge) eigenvectors.

\medskip
\noindent
Let $\mathbf g\in\mathbb R^{N}$ be the newly added off–diagonal column in the update $H_{N}\to H_{N+1}$, and define $g_i^{(N)} \;:=\; \big(\mathbf{u}_i^{(N)}\big)^{\!\top}\mathbf{g}$ for $i=1,\dots,\ell$. The one–step overlap between the edge eigenspaces at sizes $N+1$ and $N$ is $M_N := \big(U_\ell^{(N+1)}\big)^\top U_\ell^{(N)} \in \mathbb R^{\ell\times \ell}$,with consecutive spaces identified via the canonical embedding of Convention~\ref{conv:canonical-embedding}.

\medskip
\noindent
Introduce the Airy rescaling at the left edge
\[
\alpha_i^{(N)} := N^{1/6}\!\bigl(\lambda_{i}^{(N)} + 2\sqrt{N}\bigr),\qquad i=1,\dots,\ell,
\]
so that
\[
\lambda_{i}^{(N)}-\lambda_{j}^{(N)}
= N^{-1/6}\,(\alpha_i^{(N)}-\alpha_j^{(N)})
+ o_{\mathbb{P}}\!\big(N^{-1/6}\big),
\qquad 1 \le i,j \le \ell.
\]

\medskip
\noindent
Applying the ratios identity \eqref{eq:ratios} to $(H_{N},H_{N+1})$ gives
\begin{equation}\label{eq:M-ratio}
\frac{(M_N)_{ij}}{(M_N)_{ii}}
\;=\;
\frac{\Delta_{i,i}^{(N)}}{\Delta_{i,j}^{(N)}}\cdot \frac{g_j^{(N)}}{g_i^{(N)}},
\qquad
\Delta_{i,j}^{(N)} := \lambda_{i}^{(N+1)}-\lambda_{j}^{(N)},\quad 1 \le i,j \le \ell .
\end{equation}

\noindent
By adjacent-step alignment \eqref{eq:adjacent-alignment}, \((M_N)_{ii}=1+O_{\mathbb P}(N^{-2/3})\), and by the edge increment \eqref{eq:soft-edge-gap},
\begin{subequations}\label{eq:delta-entries}
\begin{align}
\Delta_{i,i}^{(N)} &=-\frac{|g_i^{(N)}|^2}{\sqrt{N}}+o_{\mathbb P}\!\big(N^{-1/2}\big),
\label{eq:delta-diag}\\
\Delta_{i,j}^{(N)} &=N^{-1/6}(\alpha_i^{(N)}-\alpha_j^{(N)})+o_{\mathbb P}\!\big(N^{-1/6}\big).
\label{eq:delta-offdiag}
\end{align}
\end{subequations}

\noindent
Combining \eqref{eq:M-ratio} and \eqref{eq:delta-entries} and, for $i\neq j$,
\[
(M_N)_{ij}
= \frac{-1}{N^{1/3}}\,
\frac{\overline{g_i^{(N)}}\,g_j^{(N)}}{\alpha_i^{(N)}-\alpha_j^{(N)}}\,\bigl(1+o_{\mathbb P}(1)\bigr),
\qquad
(M_N)_{ii}=1+O_{\mathbb P}\!\big(N^{-2/3}\big).
\]
Equivalently, by \eqref{eq:gaussian-edge-coords} we have
\((g_i^{(N)})_{i=1}^\ell=(g_i)_{i=1}^\ell+o_{\mathbb{P}}(1)\), and by the convergence to the
Airy process,
\((\alpha_i^{(N)})_{i=1}^\ell=(\alpha_i)_{i=1}^\ell+o_{\mathbb{P}}(1)\).
Consequently,
\[
M_N \;=\; I_\ell \;+\; \frac{1}{N^{1/3}}\,A_N \;+\; o_{\mathbb{P}}(N^{-1/3}),
\qquad
(A_N)_{ij} \;=\;
\begin{cases}
\dfrac{\overline{g_i}\,g_j}{\alpha_j-\alpha_i}, & i\neq j,\\[1ex]
0, & i=j.
\end{cases}
\]
so $A_N$ is skew–symmetric in the GOE case $(\beta=1)$ and skew–Hermitian in the GUE case $(\beta=2)$.

\hfill\qed

\subsubsection{Bulk Regime}
\noindent
In this section, we establish the behavior of the overlap matrix in regimes located in the spectral bulk.

\noindent
Fix \(E\in(-2,2)\) we denote \( i_N(E) \) the smallest index such that \( \lambda^{(N)}_{i_N(E)} \geq E \sqrt{N} \) and define the anchor
\[
i_N(E):=\min\bigl\{\,k\in\{1,\dots,N\}:\ \lambda_k^{(N)}\ge E\sqrt{N}\,\bigr\}.
\] 
We rescale and relabel the eigenvalues
\begin{equation}\label{eq:mu-rescale}
\mu_{j}^{(N)} \;\coloneqq\; \sqrt{4-E^2}\sqrt{N}\,\Bigl(\lambda^{(N)}_{\,j + i_N(E)} - E\sqrt{N}\Bigr),
\qquad 1 \le j + i_N(E) \le N .
\end{equation}
and relabel the weights
\begin{equation}\label{eq:def-local-g}
g_{j}^{(N)} :=
\begin{cases}
\left\langle \mathbf{u}^{(N)}_{\,j+i_N(E)},\, \mathbf{g}^{(N)} \right\rangle, & 1 \le j+i_N(E) \le N,\\
0, & \text{otherwise}.
\end{cases}
\end{equation}

\begin{proposition}[Almost-sure convergence of the input pairs]\label{prop:as-convergence}
For each fixed $j\in\mathbb{Z}$, the random pairs \eqref{eq:mu-rescale},\eqref{eq:def-local-g} $(\mu_{j}^{(N)}, g_{j}^{(N)})$ converge almost surely, as $N\to\infty$, to a limiting pair $(\mu_{j}^{(s=0)}, g_{j}^{(s=0)})$; that is,
\[
(\mu_{j}^{(N)}, g_{j}^{(N)}) \xrightarrow{\mathrm{a.s.}} (\mu_{j}^{(s=0)}, g_{j}^{(s=0)}) \qquad (N\to\infty).
\]
Moreover, the limiting points $(\mu_{j}^{(s=0)})_{j\in\mathbb{Z}}$ form the sine--$\beta$ point process, and $(g_{j}^{(s=0)})_{j\in\mathbb{Z}}$ are independent Gaussian random variables, independent of $(\mu_{j}^{(s=0)})_{j\in\mathbb{Z}}$.
\end{proposition}
\noindent
This assumption is motivated by arguments analogous to those in~\cite{huang2022}.

\bigskip
\noindent
Recall that the eigenvalues of $H_{N+1}$ are precisely the zeros of the secular equation~\eqref{eq:secular}. For $1\le j\le N$, the overlaps admit the representation
\begin{equation}\label{eq:overlap-gaps}
\Omega^{(N)}_{i,j}
= \frac{1}{\sqrt{|f'\!\bigl(\lambda^{(N+1)}_i|\bigr)}}\,
   \frac{g_j^{(N)}}{\Delta_{i,j}^{(N)}},
\qquad
\Delta_{i,j}^{(N)} := \lambda_{i}^{(N+1)}-\lambda_{j}^{(N)}.
\end{equation}
Up to the overall rescaling factor $\bigl(|f'\!\bigl(\lambda^{(N+1)}_i\bigr)\bigr)^{-1/2}$, the overlaps are a ratio of Gaussian coordinates and of eigenvalue gaps; see Lemma ~\eqref{lem:ratios}. Under Proposition~\ref{prop:as-convergence}, this ratio converges almost surely as $N\to\infty$ (for each fixed $i,j$).

\noindent
The remaining normalization factor is
\[
\upsilon_{i}^{(N)} \coloneqq \bigl(|f'\!\bigl(\lambda^{(N+1)}_i|\bigr)\bigr)^{-1/2}.
\]
We now seek the limit of $\upsilon_{i}^{(N)}$ (as $N\to\infty$) under the hypothesis of Proposition~\ref{prop:as-convergence}. Our goal is to identify an \(N\)-independent limiting expression that captures the relevant scaling.  We investigate the asymptotic behavior, as \(N \to \infty\), of the quantity

\begin{equation}\label{eq:diag-overlap-factor}
\bigl(\Delta_{i,i}^{(N)}\bigr)^{2}\, \bigl(\upsilon_{i}^{(N)}\bigr)^{-2}
= \bigl(\Delta_{i,i}^{(N)}\bigr)^{2}
\left[\,1+\sum_{j=1}^{N} \frac{\lvert g_{j}^{(N)}\rvert^{2}}{\bigl(\Delta_{i,j}^{(N)}\bigr)^{2}}\,\right].
\end{equation}
As noted in \eqref{eq:diag-overlap}, this quantity is of the same order of magnitude as the diagonal overlap, which is \emph{a priori} \(O(1)\).

\noindent
The pathological term in \eqref{eq:diag-overlap-factor} corresponding to the projection onto the new direction, namely \((\Delta_{i,i}^{(N)})^{2}\), satisfies
\((\Delta_{i,i}^{(N)})^{2}\le (\lambda_{i+1}^{(N)}-\lambda_{i}^{(N)})^{2}\) (by the Cauchy interlacing theorem).
Since the rescaling factor is of order \(N^{-1}\), this contribution becomes negligible as \(N\to\infty\), almost surely. Under the rescaling \eqref{eq:mu-rescale}, we now aim to quantify the limit of

\begin{equation}\label{eq:rescaled-norm-sum}
\bigl(\mu_{i}^{(N+1)}-\mu_{i}^{(N)}\bigr)^{2}
\sum_{j\in\mathbb{Z}}
\left(\frac{|g_{j}^{(N)}|}{\mu_{i}^{(N+1)}-\mu_{j}^{(N)}}\right)^{2}.
\end{equation}

\begin{remark}
Unlike in the eigenvalue case, where the sum to be estimated is taken in the principal value sense, the summation in our setting diverges in the bulk. The prefactor is therefore essential to properly capture the limiting behavior of this quantity. The following proposition (Proposition~\ref{prop:rescaled}) justifies interchanging the limit and the sum.
\end{remark}

\begin{proposition}\label{prop:rescaled}
Fix \(E\in(-2,2)\) and, for each \(N\), let \(i_N(E):=\min\{k:\,\lambda_k^{(N)}\ge E\sqrt{N}\}\).
For fixed offsets \(u\in\mathbb{Z}\), set \(i:=i_N(E)+u\).
Assume that the random family \(\bigl(\mu_{j}^{(N)},\,g_{j}^{(N)}\bigr)_{j\in\mathbb{Z}}\) converges almost surely, as \(N\to\infty\), to \(\bigl(\mu_{j}^{(0)},\,g_{j}^{(0)}\bigr)_{j\in\mathbb{Z}}\).
Then, almost surely,
\begin{equation}\label{eq:rescaled}
\lim_{N\to\infty}
\left(\frac{\Delta_{i,i}^{(N)}}{\upsilon_{i}^{(N)}}\right)^{2}
=
\bigl(\mu_i^{(1)}-\mu_i^{(0)}\bigr)^{2}
\sum_{j\in\mathbb{Z}}
\left(\frac{|g_j^{(0)}|}{\mu_i^{(1)}-\mu_j^{(0)}}\right)^{2},
\end{equation}
By Proposition~\ref{prop:as-convergence}, the variables \(\{g_j^{(0)}\}_{j\in\mathbb{Z}}\) are independent Gaussian random variables (independent of \(\{\mu_j^{(0)}\}_{j\in\mathbb{Z}}\)). \(\mu_i^{(1)}\) denotes the \(i\)-th solution \(z\) of the secular equation
\(
S_{\mu^{(0)},\,(g^{(0)})^{2}}(z)\;=\;-\frac{E}{2\sqrt{4-E^{2}}}\,.
\)
\end{proposition}

\begin{remark}
We observe that the right-hand side corresponds to the derivative of the function \eqref{eq:secular-S}
\(
S_{\mu,g^2}
\),
which arises in the eigenvalue process, evaluated at \(z=\mu_i^{(s+1)}\).
Since \(S_{\mu,g^2}\) is a bijection onto \(\mathbb{R}\) on each interval \([\mu_i,\mu_{i+1}]\),
define
\[
S_{\mu,g^2}(z) \coloneqq \sum_{j\in\mathbb{Z}} \frac{|g_j|^2}{\mu_j - z},
\qquad
h \coloneqq -\,\frac{E}{2\sqrt{4 - E^2}}.
\]
Then the limiting rescaled normalization factor is
\begin{equation}
\begin{aligned}
\left(\sum_{j\in\mathbb{Z}} \frac{|g_j^{(0)}|^{2}}{\bigl(\mu_i^{(1)}-\mu_j^{(0)}\bigr)^{2}}\right)^{-1}
&= \frac{1}{\,S'_{\mu,\,g^{2}}\!\left(S_{\mu,\,g^{2}}^{-1}(h)\right)_i}
&= \left(S_{\mu,\,g^{2}}^{-1}\right)^{\prime}(h)_i.
\end{aligned}
\end{equation}

\end{remark}

\paragraph*{Proof of Theorem~\ref{thm:bulk-case}}\mbox{}\\
\noindent
Denote by \(\mathbf{u}_i^{(N)}\) and \(\mathbf{u}_i^{(N+1)}\) the normalized eigenvectors associated with
\(\lambda_i^{(N)}\) and \(\lambda_i^{(N+1)}\) of the Wigner matrices \(H_{N}\) and \(H_{N+1}\), respectively.
Let \(\mathbf{g}^{(N)}\in\mathbb{R}^{N}\) be the vector formed by the last added column (excluding the diagonal entry) in the passage from \(H_{N}\) to \(H_{N+1}\).

\medskip
\noindent
Fix \(E\in(-2,2)\) and define the anchor
\[
i_N(E):=\min\bigl\{\,k\in\{1,\dots,N\}:\ \lambda_k^{(N)}\ge E\sqrt{N}\,\bigr\}.
\]
We index locally around the anchor \(i_N(E)\) by setting \(i:=i_N(E)+u\) and \(j:=i_N(E)+v\), where \(u,v\in\mathbb{Z}\) are fixed integer offsets relative to the anchor. Then, using \eqref{eq:master}, the overlap takes the form, where we set
\(g_i^{(N)}:=\mathbf{g}^{(N)\!\top}\mathbf{u}_i^{(N)}\) (the \(i\)-th coordinate of the newly added row/column when passing from \(H_N\) to \(H_{N+1}\), expressed in the eigenbasis \(\{\mathbf{u}_i^{(N)}\}\)),
and
\(\Delta_{i,j}^{(N)}:=\lambda_i^{(N+1)}-\lambda_j^{(N)}\) (the inter-step eigenvalue gap between levels \(N{+}1\) and \(N\)):

\begin{equation}
\Omega_{ij}^{(N)} =
\frac{1}{
    \bigl|\Delta_{i,i}^{(N)}\bigr|\,
    \left(
        1 + \sum_{k \in \mathbb{Z}}
        \left( \frac{ g_k^{(N)} }{ \Delta_{i,k}^{(N)} } \right)^{\!2}
    \right)^{\!1/2}
}
\;\cdot\;
\frac{ g_j^{(N)} \,\bigl|\Delta_{i,i}^{(N)}\bigr| }{ \Delta_{i,j}^{(N)} }
\label{eq:overlap-master-factorized-Delta}
\end{equation}

\noindent
Under Proposition~\ref{prop:as-convergence}, we pass to the rescaled local coordinate defined by $\mu_{j}^{(N)}:=\sqrt{4-E^2}\sqrt{N}\,(\lambda^{(N)}_{\,j+i_N(E)}-E\sqrt{N})$ for $1\le j+i_N(E)\le N$, which admits an almost-sure limit as $N\to\infty$. With this change of variables, the overlap takes the factorized form
\begin{equation}
\Omega_{ij}^{(N)} =
\frac{1}{
    \bigl|\Delta_{i,i}^{(N)}\bigr|\,
    \left(
        1 + \sum_{k \in \mathbb{Z}}
        \left( \frac{ g_k^{(N)} }{ \Delta_{i,k}^{(N)} } \right)^{\!2}
    \right)^{\!1/2}
}
\;\cdot\;
\frac{ g_j^{(N)} \,\bigl| \mu_i^{(N+1)} - \mu_i^{(N)} \bigr| }{ \mu_i^{(N+1)} - \mu_j^{(N)} } ,
\label{eq:overlap-master-factorized-mu}
\end{equation}

\noindent
Under Proposition~\ref{prop:as-convergence}, the right-hand side converges \emph{almost surely} to its limiting counterpart by continuity of finite algebraic operations. The behavior of the rescaling factor is handled in Proposition~\ref{prop:rescaled}. Consequently, standard limit arguments yield the \emph{almost sure} limit of the overlaps, confirming the claimed formula:
\[
\big\langle \mathbf{u}_i^{(N+1)},\, \mathbf{u}_j^{(N)} \big\rangle
\;\xrightarrow[N\to\infty]{\mathrm{a.s.}}\;
\left( \sum_{j \in \mathbb{Z}} 
       \left( \frac{|g_j^{(0)}|}{\mu_i^{(1)} - \mu_j^{(0)}} \right)^{2} 
\right)^{-1/2} \,\frac{g_j^{(0)}}{\mu_j^{(1)} - \mu_i^{(0)}}
\]

\noindent
The theorem follows from the preceding remark together with the almost-sure limit established above.
\hfill\qed

\begin{remark}

\noindent
Proceeding as in the proof of Theorem~\ref{thm:bulk-case} and adopting the same notation, we also obtain joint convergence along \emph{paired} indices. Fix a finite family of pairs \(\mathcal{P}=\{(i_1,j_1),\dots,(i_r,j_r)\}\subset\mathbb{Z}^2\). Then
\begin{equation}\label{eq:joint-overlaps-limit-pairs}
\bigl(\Omega_{i_N(E)+i_u, i_N(E)+j_u}^{(N)}\bigr)_{1\le u\le r}
\;\xrightarrow[N\to\infty]{\mathrm{a.s.}}\;
\left(
\left[\sum_{j\in\mathbb{Z}}
\left(\frac{|g_j^{(0)}|}{\mu_{i_u}^{(1)}-\mu_j^{(0)}}\right)^2\right]^{-1/2}
\frac{g_{j_u}^{(0)}}{\mu_{j_u}^{(1)}-\mu_{i_u}^{(0)}}
\right)_{1\le u\le r}.
\end{equation}
\noindent
Here \(\bigl(\mu_i^{(0)}\bigr)_{i\in\mathbb{Z}}\) is the \(\mathrm{Sine}_\beta\) point process, \(\bigl(g_i^{(0)}\bigr)_{i\in\mathbb{Z}}\) are i.i.d.\ \(\mathcal{N}_\beta(0,1)\) and independent of \(\bigl(\mu_i^{(0)}\bigr)\), \(\bigl(\mu_i^{(1)}\bigr)_{i\in\mathbb{Z}}\) are the ordered solutions \(z\) to the secular equation \(S_{\mu,g^2}(z)=h\), and \(v_u^{(0)}=\sqrt{\bigl(S_{\mu,g^2}^{-1}\bigr)'(h)_u}\) (with \(h\) as in \cite{huang2022}) using (Convention~\ref{conv:inverse-stieltjes-branches}).

\end{remark}

\paragraph*{Proof of Proposition~\eqref{prop:rescaled}}\mbox{}\\

\noindent
Let $z\in\mathbb{R}$. We carry out all computations for a real parameter $z$, and only at the end do we substitute the (real) solution of the secular equation from~\cite{huang2022}. We split the target sum into four parts: (i) a finite block that converges almost surely; (ii) a sum–integral comparison for the remaining terms; (iii) the integral of the remaining tail; and (iv) a comparison in which the random variables $|g_j^{(N)}|^2$ are replaced by their expectations. The last three contributions are uniformly bounded in $k$. Letting $N\to\infty$ first, the sum decomposes into the limit of the finite block plus a remainder depending only on $k$, and this remainder tends to $0$ as $k\to\infty$, which determines the limit. We start from equation~\eqref{eq:rescaled-norm-sum} with \(i:=i_N(E)\) (the computation for \(i=i_N(E)+u\) is analogous), using the notation introduced in \eqref{eq:mu-rescale} and \eqref{eq:def-local-g}.

\begin{align*}
\sum_{j \in \mathbb{Z}} \bigl|g_{j}^{(N)}\bigr|^2 \frac{(\mu_{0}^{(N)} - z)^2}{(\mu_{j}^{(N)} - z)^2} 
&=
\underbrace{\sum_{|j| \leq k} \bigl|g_{j}^{(N)}\bigr|^2 \frac{(\mu_{0}^{(N)} - z)^2}{(\mu_{j}^{(N)} - z)^2}}_{(1)}
+
\underbrace{(\lambda_{i_N(E)}^{(N)} - z)^2 \int_{|x - E| \geq \tfrac{2\pi k}{N\sqrt{4-E^2}}} \frac{\rho_{sc}(x)}{(x - z)^2} \, dx}_{(2)}
\\
&\quad+
\underbrace{\left(
\sum_{|j| \geq k} \frac{(\mu_{0}^{(N)} - z)^2}{(\mu_{j}^{(N)} - z)^2}
- 
(\lambda_{i_N(E)}^{(N)} - z)^2 \int_{|x - E| \geq \tfrac{2\pi k}{N\sqrt{4-E^2}}} \frac{\rho_{sc}(x)}{(x - z)^2} \, dx
\right)}_{(3)}
\\
&\quad+
\underbrace{\sum_{|j| \geq k} \frac{(\mu_{0}^{(N)} - z)^2}{(\mu_{j}^{(N)} - z)^2} \bigl( \,|g_{j}^{(N)}|^2 - 1 \,\bigr)}_{(4)}.
\end{align*}

\begin{itemize}

\item[1.] 
We begin by analyzing the finite summation over indices \(|j| \leq k\). Since the limit of a finite sum equals the sum of the limits, we obtain the convergence:
\[
\lim_{N \to \infty} \sum_{|j| \leq k} \bigl|g_{j}^{(N)}\bigr|^2 \frac{(\mu_{0}^{(N)} - z)^2}{(\mu_{j}^{(N)} - z)^2}
= \sum_{|j| \leq k} |g_i^{(0)}|^2 \frac{(\mu_0^{(0)} - z)^2}{(\mu_j^{(0)} - z)^2}.
\]

\item[2.]
We now estimate the integral term that appears in the bulk regime. Using the asymptotics near the simple pole at $x=z$ and the spacing estimate from \cite{huang2022},
\(
\mu_{k}^{(N)}=2\pi\,k+O(k^{3/4}),
\)
we set the cutoff $\varepsilon_{N,k} = \tfrac{2\pi k}{N\sqrt{4-E^2}}$ and write $y=x-z$. Then
\[
\int_{|x-z|\ge \tfrac{2\pi k}{N}} \frac{\rho_{\mathrm{sc}}(x)}{(x-z)^{2}}\,dx
\;\lesssim\; \rho_{\mathrm{sc}}(z)\!\int_{|y|\ge \varepsilon_{N,k}} \frac{dy}{y^{2}}
\;=\; \rho_{\mathrm{sc}}(z)\,\frac{2}{\varepsilon_{N,k}}
\;\lesssim\; \rho_{\mathrm{sc}}(z)\,\frac{N}{k}.
\]
(Here we used the elementary bound $\displaystyle \int_{|y|\ge\varepsilon} y^{-2}\,dy=2/\varepsilon\lesssim 1/\varepsilon$ as $\varepsilon\to 0$.) Hence
\[
(\lambda_{i_N(E)}^{(N)}-z)^{2}\!
\int_{|x-z|\ge \tfrac{2\pi k}{N\sqrt{4-E^2}}} \frac{\rho_{\mathrm{sc}}(x)}{(x-z)^{2}}\,dx
\;\lesssim\;
(\lambda_{i_N(E)}^{(N)}-z)^{2}\cdot \frac{\rho_{\mathrm{sc}}(z)\,N}{2\pi k}.
\]
Recalling that
\(
(\lambda_{i_N(E)}^{(N)}-z)^{2}
= N^{-1}(\mu_{0}^{(N)}-z)^{2}\lesssim N^{-1},
\)
we conclude the bound
\[
(\lambda_{i_N(E)}^{(N)}-z)^{2}
\int_{|x-z|\ge \tfrac{2\pi k}{N\sqrt{4-E^2}}} \frac{\rho_{\mathrm{sc}}(x)}{(x-z)^{2}}\,dx
\;\lesssim\; k^{-1}.
\]

\begin{remark}
In particular, the sum of squared overlaps outside the $k$-diagonal band decays like $k^{-1}$, uniformly in $N$:
\[
\sum_{|i-j|\ge k} \bigl|\Omega_{ij}^{(N)}\bigr|^{2} \;\lesssim\; k^{-1}.
\]
\end{remark}

\item[3.]
We next show that the difference between the discrete sum and its continuous (semicircle) integral becomes negligible as \(k\to\infty\).
Using Stieltjes integration by parts (summation by parts) and comparing the empirical spectral measure with its limiting semicircle law (as in \cite{huang2022}), we obtain
\[
\left|
\sum_{j\ge k}\frac{1}{(\mu_j^{(N)}-z)^2}
-\frac{1}{N}\int_{E+\mu_k^{(N)}/N\sqrt{4-E^2}}^{\infty}\frac{\rho_{sc}(x)}{(x-z)^2}\,dx
\right|
\;\lesssim\;
\int_{\mu_k^{(N)}}^{\infty}\frac{(1+x)^{3/4}}{(x-z)^3}\,dx
\;\lesssim\; k^{-5/4}.
\]

\item[4.]
We bound the error incurred by replacing $|g_j^{(N)}|^{2}$ with its mean $1$ as $k\to\infty$.
By \cite{huang2022}, there exists $\sigma>0$ such that, uniformly in the lower cutoff and
\emph{almost surely},
\[
\sup_{\ell\ge k}
\left|
\sum_{j\ge \ell}\frac{|g_{j}^{(N)}|^2-1}{\mu_{j}^{(N)}-z}
\right|
\;\lesssim\;
k^{-(1-1/\sigma)/4}+k^{-1/\sigma},
\]

Set 
\(
a_j:=\dfrac{|g_{j}^{(N)}|^2-1}{\mu_{j}^{(N)}-z}
\)
and 
\(
b_j:=\dfrac{1}{\lvert \mu_{j}^{(N)}-z\rvert},
\)
so that $(b_j)_{j\ge k}$ is positive and decreasing in the bulk. Writing
\(A_\ell:=\sum_{j\ge \ell}a_j\),
Abel’s summation gives
\[
\sum_{j\ge k}\frac{|g_{j}^{(N)}|^2-1}{(\mu_{j}^{(N)}-z)^2}
=\sum_{j\ge k} a_j b_j
= A_k b_k + \sum_{j>k} A_j\,(b_j-b_{j-1}).
\]
Taking absolute values and using that $b_j$ decreases to $0$,
\[
\left|\sum_{j\ge k}\frac{|g_{j}^{(N)}|^2-1}{(\mu_{j}^{(N)}-z)^2}\right|
\le
\bigl(A_k + \sup_{j\ge k}|A_j|\bigr)\, b_k
\;\lesssim\;
\frac{1}{\lvert \mu_k^{(N)}-z\rvert}
\Bigl(k^{-(1-1/\sigma)/4}+k^{-1/\sigma}\Bigr).
\]
By bulk rigidity/level spacing, $\lvert \mu_k^{(N)}-z\rvert\gtrsim k$, hence $b_k\lesssim k^{-1}$, and therefore
\[
\left|\sum_{j\ge k}\frac{|g_{j}^{(N)}|^2-1}{(\mu_{j}^{(N)}-z)^2}\right|
\;\lesssim\; k^{-1}\xrightarrow[k\to\infty]{}0.
\]
Thus the replacement $|g_j^{(N)}|^2\mapsto 1$ produces a tail error that vanishes as $k\to\infty$.

\end{itemize}

\bigskip

\noindent
Combining the preceding estimates, we obtain the following limit.
Let \(z\) be the solution to the secular equation lying in the open interval
\((\mu_{i}^{(0)},\,\mu_{i+1}^{(0)})\). Then, by Theorem~1.4 of Huang \cite{huang2022}, we have
\begin{equation}
\lim_{N \to \infty} 
\left( \frac{\lambda_{i}^{(N+1)} - \lambda_{i}^{(N)}}{\upsilon_{i}^{(N)}} \right)^{2}
=
\bigl(\mu_i^{(1)} - \mu_i^{(0)}\bigr)^{2}
\left( \sum_{j \in \mathbb{Z}} 
       \left( \frac{|g_j^{(0)}|}{\mu_i^{(1)} - \mu_j^{(0)}} \right)^{2} 
\right).
\end{equation}

\qed

\subsubsection{Regimes for the remaining overlap entries}

\noindent Building on the preceding computations, we state three additional asymptotic regimes for the remaining overlap entries as immediate corollaries. The results depend on the spectral band under consideration; accordingly, we introduce the quantile parameter $\lambda(y_i)$ defined by $\int_{-\infty}^{\lambda(y_i)} \rho_{\mathrm{sc}}(x)\,dx = y_i$ for $y_i\in(0,1)$. Adopting the same notation as in the statement of Theorem~\ref{thm:bulk-case}, we obtain the following marginal laws; the joint laws can be derived analogously to~\eqref{eq:joint-overlaps-limit-pairs}. We recall Convention~\ref{conv:inverse-stieltjes-branches}.

\begin{itemize}
  \item One index at the edge and one in the bulk. For \(i=O(1)\) and \(j=yN\) with fixed \(y\in(0,1)\),
  \begin{equation}\label{eq:edge-bulk}
    N\,\Omega^{(N)}_{i,j}\ \xrightarrow{d}\ 
    \frac{\overline{g_i}\,g_j}{\lambda(0)-\lambda(y)} .
  \end{equation}
  I.e., a deterministic \(1/(\text{edge–bulk})\) profile at scale \(N^{-1}\). The same formula holds at the upper edge upon replacing \(\lambda(0)\) by the upper-edge value (namely \(2\) in semicircle units).

  \item Two bulk indices with macroscopic separation: for \(i=y_1N\), \(j=y_2N\) with distinct \(y_1,y_2\in(0,1)\),
  \begin{equation}\label{eq:bulk-off}
    N\,\Omega^{(N)}_{i,j}\ \xrightarrow{d}\ 
    g_j\,
    \frac{\sqrt{\bigl(S^{-1}_{\mu,\,g^2}\bigr)'(h)_i}}
         {\lambda(y_j)-\lambda(y_i)} ,\qquad
    h=-\frac{E}{2\sqrt{4-E^2}} .
  \end{equation}
  Thus off–diagonal bulk overlaps are of order \(N^{-1}\) with the usual \(1/(\lambda(y_j)-\lambda(y_i))\) decay (here \(S^{-1}_{\mu,\,g^2}\) is as defined in the bulk analysis; the derivative acts componentwise).

  \item For the overlap with the freshly added basis direction, for \(i=yN\) with fixed \(y\in(0,1)\),
  \begin{equation}\label{eq:new-direction}
    \sqrt{N}\,\Omega^{(N)}_{i,\,N+1}
    \ \xrightarrow{d}\ 
    \frac{\sqrt{\bigl(S^{-1}_{\mu,\,g^2}\bigr)'(h)_i}}{\lambda(y)} ,\qquad
    h=-\frac{E}{2\sqrt{4-E^2}} .
  \end{equation}
  Hence these entries are of order \(N^{-1/2}\).
\end{itemize}

\subsection{Master Equations for Eigenvalues and Eigenvectors in the Wishart case}

\begin{proposition}[Eigenvalue master equation — Wishart case]\label{prop:wishart-master}
Fix integers $N, T \ge 1$, and let $\{\lambda_i\}_{i=1}^{N}$ denote the (positive) eigenvalues of $X_{N}^{\!\top} X_{N}$.
Consider the rank-one extension
\[
X_{N+1} \;=\; \big( X_N \;\; g \big),
\qquad g \in \mathbb{R}^{T},
\]
that is, $X_{N+1}$ is obtained by appending the column vector $g$ to $X_N$.
Then the eigenvalues $z$ of $X_{N+1}^{\!\top} X_{N+1}$ are the zeros of the \emph{secular equation}
\begin{equation}\label{eq:wishart-secular-simplified}
F_N(z)
\;=\;
1
+ \sum_{i=1}^{N}
\frac{|\langle v_i, g \rangle|^2}{\lambda_i - z}
- \frac{\Gamma}{z}
\;=\;0,
\end{equation}
where ${v_i}$ are the eigenvectors associated with the nonzero eigenvalues of $X_N X_N^{\top}$, and
\(
\Gamma
:= \sum_{j=1}^{T-N} |\langle w_j, g \rangle|^2,
\)
with $\{w_j\}$ forming an orthonormal basis of the null space of $X_{N} X_{N}^{\!\top}$.
\end{proposition}

\begin{proposition}[Eigenvector master equation — Wishart case]\label{prop:wishart-master-eigenvector}
Fix integers $N, T \ge 1$, and let $\{\lambda_i^{(N)}\}_{i=1}^{N}$ 
denote the (ordered, positive) eigenvalues of the Wishart matrix 
\(
W^{(N)} := X_{N}^{\!\top} X_{N} \in \mathbb{R}^{N\times N}.
\)
Let $\{\mathbf{u}_i^{(N)}\}_{i=1}^{N}$ be the associated orthonormal eigenvectors. 
Consider the rank–one extension
\[
X_{N+1} = \big( X_N \;\; g \big),
\qquad g \in \mathbb{R}^{T},
\]
that is, $X_{N+1}$ is obtained by appending the column vector $g$ to $X_N$.
Let $W^{(N+1)} = X_{N+1}^{\!\top} X_{N+1}$, with eigenvalues 
$\{\lambda_i^{(N+1)}\}_{i=1}^{N+1}$ and eigenvectors $\{\mathbf{u}_i^{(N+1)}\}$.
We define
\(
g_i^{(N)} \;=\; \langle v_i^{(N)},\, g \rangle,
\)
where $v_i^{(N)}$ denotes the eigenvector associated with the eigenvalue 
$\lambda_i^{(N)}$ of the matrix $X_N X_N^{\!\top}$.

\smallskip
\noindent
Then the overlap coefficients
\(
\Omega^{(N)}_{ij} := 
\big\langle \mathbf{u}_i^{(N+1)},\, \mathbf{u}_j^{(N)} \big\rangle
\)
satisfy the relation
\begin{equation}\label{eq:wishart-ratio-identities-compact}
\begin{aligned}
\frac{\Omega^{(N)}_{k j}}{\Omega^{(N)}_{k i}} 
&= 
\frac{g_j^{(N)}}{g_i^{(N)}}
\sqrt{\frac{\lambda_j^{(N)}}{\lambda_i^{(N)}}}
\frac{\Delta_{k,i}^{(N)}}{\Delta_{k,j}^{(N)}},
\qquad 1 \le i,j \le N, \qquad
\frac{\Omega^{(N)}_{k,\,N+1}}{\Omega^{(N)}_{k i}} 
&= 
\frac{\Delta_{k,i}^{(N)}}{g_i^{(N)}}
\sqrt{\frac{1}{\lambda_i^{(N)}}}.
\end{aligned}
\end{equation}
where
\(
\Delta_{k,i}^{(N)} := \lambda_k^{(N+1)} - \lambda_i^{(N)}.
\)

\end{proposition}

\begin{proof}[Proof of Proposition~\ref{prop:wishart-master}]
Write $X:=X_N\in\mathbb R^{T\times N}$ and $G:=X_{N+1}=(X\;\;g)$. Then
\[
W^{(N)}=X^\top X,\qquad 
W^{(N+1)}=G^\top G=
\begin{bmatrix}
X^\top X & X^\top g\\ g^\top X & g^\top g
\end{bmatrix}.
\]
Fix $z\in\mathbb C$ with $z\notin\{\lambda_1,\dots,\lambda_N\}$ and apply the Schur–complement formula
to $W^{(N+1)}-zI$:
\[
\det\!\big(W^{(N+1)}-zI\big)
= \det\!\big(W^{(N)}-zI\big)\;
\det\!\Big(g^\top g - z - g^\top X\,(X^\top X-zI)^{-1}X^\top g\Big).
\]
Hence $z$ is an eigenvalue of $W^{(N+1)}$ iff the scalar factor
\[
\Phi_N(z):=g^\top g - z - g^\top X\,(X^\top X-zI)^{-1}X^\top g
\]
vanishes. Let $\{v_i\}_{i=1}^N$ be orthonormal eigenvectors of $X X^\top$ with $X X^\top v_i=\lambda_i v_i$,
and let $\{w_j\}_{j=1}^{T-N}$ be an orthonormal basis of $\ker(X^\top)$ (equivalently, of the zero eigenspace of
$X X^\top$). Decompose
\(
g=\sum_{i=1}^N \langle v_i,g\rangle v_i + \sum_{j=1}^{T-N}\langle w_j,g\rangle w_j.
\)
Since $X(X^\top X - zI)^{-1}X^\top=\sum_{i=1}^N \frac{\lambda_i}{\lambda_i-z}\, v_i v_i^\top$,
we compute
\[
g^\top X\,(X^\top X-zI)^{-1}X^\top g
=\sum_{i=1}^N \frac{\lambda_i}{\lambda_i-z}\,|\langle v_i,g\rangle|^2,
\qquad
g^\top g=\sum_{i=1}^N |\langle v_i,g\rangle|^2+\sum_{j=1}^{T-N} |\langle w_j,g\rangle|^2.
\]
Therefore
\[
\Phi_N(z)
=-z\left( 1 + \sum_{i=1}^N \frac{|\langle v_i,g\rangle|^2}{\lambda_i-z}
-\frac{\Gamma}{z}\right),
\qquad
\Gamma:=\sum_{j=1}^{T-N} |\langle w_j,g\rangle|^2.
\]
Since $-z\neq 0$ for any nonzero eigenvalue, the eigenvalue condition $\Phi_N(z)=0$ is equivalent to
\[
F_N(z):=1+\sum_{i=1}^N \frac{|\langle v_i,g\rangle|^2}{\lambda_i-z}-\frac{\Gamma}{z}=0,
\]
which is exactly \eqref{eq:wishart-secular-simplified}. This proves the claim.
\end{proof}

\begin{proof}[Proof of Proposition~\ref{prop:wishart-master-eigenvector}]
Work in the orthonormal basis $\{u_1^{(N)},\dots,u_N^{(N)}\}$ of eigenvectors of $W^{(N)}=X^\top X$
and complete it by the $(N+1)$-st basis vector corresponding to the new column.
As recalled in the statement (see also the construction preceding the proposition), in this basis
$W^{(N+1)}$ has the arrowhead form
\[
W^{(N+1)}=
\begin{bmatrix}
\lambda_1^{(N)} & & & 0 & \sqrt{\lambda_1^{(N)}}\,g_1^{(N)}\\
& \lambda_2^{(N)} & & & \sqrt{\lambda_2^{(N)}}\,g_2^{(N)}\\
& & \ddots & & \vdots\\
0 & & & \lambda_N^{(N)} & \sqrt{\lambda_N^{(N)}}\,g_N^{(N)}\\
\sqrt{\lambda_1^{(N)}}\,g_1^{(N)} & \sqrt{\lambda_2^{(N)}}\,g_2^{(N)} & \cdots &
\sqrt{\lambda_N^{(N)}}\,g_N^{(N)} & \;g^\top g
\end{bmatrix},
\]
where $g_i^{(N)}=\langle v_i^{(N)},g\rangle$ and $v_i^{(N)}$ is the eigenvector of $X_N X_N^\top$
associated with $\lambda_i^{(N)}$.

\noindent
This completes the proof, upon applying Lemma~\ref{lem:ratios} to the matrix $W^{(N+1)}$.

\end{proof}

\subsection{Asymptotic Behavior of Masters Equations in the Wishart case}

\noindent
The bulk and edge universality for Wishart matrices have been established in a series of works, beginning with Ben Arous and Péché \cite{benarous2005}, and extended to general covariance ensembles by Tao and Vu \cite{tao2012} and Pillai and Yin \cite{pillai2014}.
The asymptotic behavior of the largest eigenvalue at the soft edge (Tracy–Widom law) was first proved by Johnstone \cite{johnstone2001} in the Gaussian case, and later generalized by Soshnikov \cite{soshnikov2002} and Bao et al. \cite{bao2015} to non-Gaussian models.

\bigskip
\noindent
Fix any $E \in \bigl((1-\sqrt{q})^2,\,(1+\sqrt{q})^2\bigr)$, 
and denote by $j^{(s)}$ the smallest index such that 
$\lambda_{j^{(s)}}^{(s)} \ge E T$. 
We rescale and relabel the eigenvalues as
\[
\mu_j^{(s)} 
= 
\lambda_{j+j^{(s)}}^{(s)} - E T,
\qquad 
1 \le j + j^{(s)} \le N + s,
\]
and relabel the weights as
\[
\gamma_j^{(s)}
= 
\bigl|\langle \mathbf{v}_{j+j^{(s)}}^{(s)},\, \mathbf{g}^{(s+1)} \rangle\bigr|^2,
\qquad 
j + j^{(s)} \in [1, N+s],
\quad 
g_j^{(s)} = 0, 
\quad 
j + j^{(s)} \notin [1, N+s],
\]
where
\[
\mathbf{g}^{(s+1)} 
= 
\bigl(g_{1,\,N+s+1},\, g_{2,\,N+s+1},\, \dots,\, g_{T,\,N+s+1}\bigr)^\top
\]
is independent of $X_{N+K}^{(s)}$.

\medskip
\noindent
We also define the quantity
\[
\Gamma_N 
= 
\sum_{i=1}^{T-N-s} 
\bigl|\langle \mathbf{w}_{i}^{(s)},\, \mathbf{g}^{(s+1)} \rangle\bigr|^2,
\]
where 
$\{\mathbf{w}_i^{(s)}\}_{i=1}^{T-N-s}$ 
forms an orthonormal basis of the null space of 
$X_N X_N^{\!\top}$.
The term $\Gamma_N$ therefore represents the squared norm of the projection of 
$\mathbf{g}^{(s+1)}$ onto the unexplored subspace spanned by the eigenvectors 
associated with the non-zero eigenvalues of $X_N X_N^{\!\top}$.

\medskip
\noindent
To complete the proof of Theorem~\ref{thm:bulk}, it remains to establish the following result.
\begin{remark}\label{ass:huang-mp}
We assume that the random objects see equations~(4.10)--(4.12) satisfy, in the Marchenko--Pastur case, the same almost sure convergence properties established by \cite{huang2022} in the Wigner case.
\end{remark}

\begin{proposition}\label{prop:limit-wishart}

We have almost surely
\[
\lim_{N \to \infty}
\left(
1 
+ \sum_{j \in \mathbb{Z}} \frac{\gamma_j^{N}}{\mu_j^{N} - z}
- \frac{\Gamma_N}{z}
\right)
=
 \sum_{j \in \mathbb{Z}} \frac{\gamma_j}{\mu_j - z}
+ \frac{E+q-1}{2E},
\]
for all $z \in \mathbb{R} \setminus \{\mu_j\}_{j \in \mathbb{Z}}$.
\end{proposition}

\bigskip
\noindent
In the \emph{hard-edge regime} corresponding to $q = 1$ and $\alpha = T - N > 1$, 
one should note that the term $\Gamma_N$ no longer converges to its mean by the central limit theorem. 
After rescaling the eigenvalues as $\mu_j = 4N \lambda_j$, we obtain the following result.

\begin{proposition}[Hard-Edge Limit]\label{prop:limit-wishart-hard}
We have, almost surely,
\[
\lim_{N \to \infty}
\left(
\frac{1}{4N}
+ 
\sum_{j \in \mathbb{Z}} \frac{\gamma_j^{N}}{\mu_j^{N} - z}
- 
\frac{\Gamma_N}{z}
\right)
=
\sum_{j \in \mathbb{Z}} \frac{\gamma_j}{\mu_j - z}
- 
\frac{\chi_{\alpha}}{z},
\]
for all $z \in \mathbb{R} \setminus \{\mu_j\}_{j \in \mathbb{Z}}$.
\end{proposition}

\noindent
The results on eigenvector overlaps follow as corollaries of the preceding propositions. In the edge regime, by using the closed-form expression of the diagonal overlap, one recovers the result of Bao~\cite{bao2025} for the Wishart ensemble. All other statements follow from analogous arguments to those developed in the case of Wigner matrices.

\end{document}